\documentclass[sn-mathphys,Numbered]{sn-jnl}

\usepackage{graphicx}%
\usepackage{multirow}%
\usepackage{amsmath,amssymb,amsfonts}%
\usepackage{amsthm}%
\usepackage{mathrsfs}%
\usepackage[title]{appendix}%
\usepackage{xcolor}%
\usepackage{textcomp}%
\usepackage{manyfoot}%
\usepackage{booktabs}%
\usepackage{algorithm}%
\usepackage{algorithmicx}%
\usepackage{algpseudocode}%
\usepackage{listings}%
\usepackage{natbib}
\usepackage{framed}
\usepackage{latexsym}

\usepackage{url}
\definecolor{newcolor}{rgb}{.8,.349,.1}
\usepackage{mathtools}
\usepackage{enumitem} 
\usepackage{ulem,soul} 
\usepackage{graphics, subcaption}
\usepackage{hyperref}
\usepackage{todonotes}
\usepackage{booktabs,siunitx} 

\newcommand{\wt}{\widetilde}

\newcommand{\R}{\mathbb{R}}

\newcommand{\norm}[1]{\left\lVert {#1} \right\rVert}

\newcommand{\D}{\Delta}

\newcommand{\aletext}[1]{\textit{\textcolor{red}{#1}}}
\newcommand{\agntext}[1]{\textit{\textcolor{cyan}{#1}}}
\newcommand\ale[2][m]{\ifx i#1 \todo[linecolor=red,backgroundcolor=red!25,bordercolor=red,inline]{#2}\else
\todo[linecolor=red,backgroundcolor=red!25,bordercolor=red]{#2}\fi}
\newcommand\agn[2][m]{\ifx i#1 \todo[linecolor=cyan,backgroundcolor=cyan!25,bordercolor=cyan,inline]{#2}\else
\todo[linecolor=cyan,backgroundcolor=cyan!25,bordercolor=cyan]{#2}\fi}

\theoremstyle{thmstyleone}%
\theoremstyle{thmstyletwo}%
\newtheorem{remark}{Remark}%

\theoremstyle{thmstylethree}%

\begin{document}

\title[A POD approach to identify and control PDEs through SDRE]{A POD approach to identify and control PDEs online  through State Dependent Riccati equations\footnotesize{ Dedicated to Maurizio: beyond a math guide on our paths}}
\author*[1]{\fnm{Alessandro} \sur{Alla}}\email{alessandro.alla@unive.it}
\author[2]{\fnm{Agnese} \sur{Pacifico}}\email{agnese.pacifico@uniroma.it}

\affil*[1]{\orgdiv{Dipartimento di Scienze Molecolari e Nanosistemi}, \orgname{Universit\`a Ca' Foscari, Venezia}, \orgaddress{\street{Via Torino, 155}, \city{Mestre}, \postcode{30170}, \country{Italy}}}

\affil*[2]{\orgdiv{Dipartimento di Matematica}, \orgname{Sapienza Universit\`a di Roma}, \orgaddress{\street{Piazzale Aldo Moro, 5}, \city{Roma}, \postcode{00185},  \country{Italy}}}


\abstract{

We address the control of Partial Differential equations (PDEs) with unknown parameters. Our objective is to devise an efficient algorithm capable of both identifying and controlling the unknown system. We assume that the desired PDE is observable provided a control input and an initial condition.  The method works as follows, given an estimated parameter configuration, we compute the corresponding control using the State-Dependent Riccati Equation (SDRE) approach. Subsequently, after computing the control, we observe the trajectory and estimate a new parameter configuration using Bayesian Linear Regression method. This process iterates until reaching the final time, incorporating a defined stopping criterion for updating the parameter configuration. We also focus on the computational cost of the algorithm, since we deal with high dimensional systems. To enhance the efficiency of the method, indeed, we employ model order reduction through the Proper Orthogonal Decomposition (POD) method. The considered problem's dimension is notably large, and POD provides impressive speedups. Further, a detailed description on the coupling between POD and SDRE is also provided. Finally, numerical examples will show the accurateness of our method across two test cases.

}


\maketitle












 
\section{Introduction}
\begin{figure}[htbp]
\centering
\begin{tikzpicture}[scale=3]
\draw (-0.25,0.44) node[centered ] {$u(t)$};
\draw (-0.25,0.17) node[centered ] {$x_0$};
 \draw [-stealth](-0.5,0.12) -- (-0.05,0.12);
  \draw [-stealth](-0.5,0.38) -- (-0.05,0.38);
\draw (0,0) -- (1,0) -- (1,0.5) -- (0,0.5) -- (0,0); 
\fill[orange](0,0) -- (1,0) -- (1,0.5) -- (0,0.5) -- (0,0); 
\draw [-stealth](1.05,0.25) -- (2,0.25);
\draw (1.51,0.31) node[centered ] {$x(t;u(t),\mu^*)$};
\draw (0.5,0.25) node[centered ] {$\mu^*$};
\end{tikzpicture}
\caption{System observation can be seen as a black box: given a control $u(t)$ and an initial state $x_0$, we can observe the trajectory $x(t;u(t),\mu^*)$ obtained with the provided inputs. Here $\mu^*$ is an unknown parameter of the PDE. System observations allow us to see the system evolution even if $\mu^*$ is unknown. With these observations we will approximate the unknown parameter. In principle any input $u$ could be used for observing the system, but we will look for a $u$ that minimizes a given cost functional.}
\label{fig:blackbox}
\end{figure}
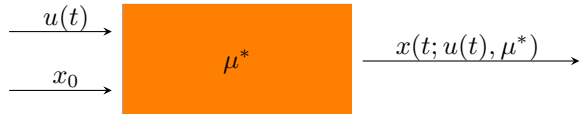

We consider the control and identification of large dimensional problems using feedback control strategies. Furthermore, we develop an efficient method by means of model order reduction. In Figure \ref{fig:blackbox}, we illustrate the configuration under investigation. For an initial state $x_0$ and a prescribed input $u(t)$ the black box, equipped with knowledge of the exact parameter configuration $\mu^*$, yields the trajectory $x(t;u(t),\mu^*)$ corresponding to the provided inputs.

Our objective is, therefore, twofold: firstly, to identify a control strategy that stabilizes the trajectory according to a specified cost functional, and secondly, to estimate the associated parameter configuration $\mu^*$. At each iteration, we compute an estimate $\wt \mu$ for the parameter configuration aiming at the convergence $\wt\mu\rightarrow \mu^*$. Later, we compute the corresponding control input $u(t; \wt \mu)$. Crucially, the estimation of $\mu^*$ relies solely on the observation of the trajectory $x(t; u(t; \wt \mu),\mu^*)$, and the entire process is designed for online execution without the necessity of computing long-term trajectories for each given input. 
This methodology, initially introduced in \cite{pacifico2021LSSC} for linear ODEs and later in \cite{RLPDE} for one dimensional PDEs, is extended to more challenging problems in this manuscript. Our current focus involves the exploration of two dimensional PDEs, and, to keep the method efficient, model reduction techniques will play a pivotal role in achieving our objective.\\ 

Indeed, optimal control problems arising from partial differential equations (PDEs) often involve the discretization of these PDEs, leading to large-dimensional systems of ordinary differential equations (ODEs). 
The large-scale nature of these problems makes traditional feedback control methods computationally intractable. The use of the Algebraic Riccati equation (ARE) in this context represents an effective means of control synthesis, offering insights into optimal strategies without an exhaustive exploration of the entire state space. As a result, the application of the ARE becomes indispensable in managing the intricacies of optimal control for large-dimensional systems arising from PDE discretization. A lot of effort has been addressed in the last two decades to approximate ARE for large scale settings. We refer to e.g. \cite{Benner.Saak.survey13,B08,Simoncini.survey.16,SSM14,S16,Binietal.book.12}.

State-Dependent Riccati Equations (SDRE, \cite{BTB00,BLT07}) present a valuable extension to the conventional Riccati framework, especially when addressing optimal control challenges for nonlinear problems. Unlike the traditional LQR approach \cite{anderson1999LQR}, SDRE incorporates state-dependent matrices, providing a more flexible approach for handling complex nonlinear dynamics. The SDRE methodology offers a tailored solution by allowing dynamic variations in the Riccati matrices based on the system's current state. This adaptability proves crucial in optimizing control strategies for diverse and intricate state-dependent behaviors. We refer to \cite{BH18} for a theoretical study on the stabilization hypothesis for SDRE in a small-scale scenario. Clearly, the computational cost of this method for large dimensional problems increases exponentially. This is due to the fact that a large scale ARE has to be computed at each time iteration. In \cite{AKS23}, the authors have addressed this issue providing numerical examples for large scale problems. Another class of methods to approximate high dimensional ARE is based on model order reduction. The goal is to project the original problem onto low dimensional subspaces. One of the first approaches goes back to \cite{AK11} where the authors have used a method based on the Proper Orthogonal Decomposition (POD, \cite{BHL93,S87}) for linear PDEs. The POD method is based on snapshots of the solution, which is the state of the art of the method. In \cite{AK11}, the authors generated snapshots and POD basis from the uncontrolled problem. More sophisticated methods have been investigated later. We refer to e.g. \cite{BS16} where the snapshots were computed from the adjoint equation and to \cite{AS18} for a comparison of several reduced approaches. A different point of view was proposed in \cite{SH18} for parametrized problems. There, the authors suggested to use a basis computed directly from the low rank solution of ARE using a greedy strategy based on certified a-posteriori error estimators.

In the context of the SDRE, a POD approach has been investigated in \cite{HW23} for quadratic nonlinear terms. In this work, we extend the study to generic nonlinear problems using also an hyper reduction based on the Discrete Empirical Interpolation Method (DEIM, \cite{CS10}) which is critical to set a reduced model independent from the original dimension of the problem. Specifically, since at each iteration one has to solve an ARE which depends on the current state space nonlinearly, the use of DEIM is decisive to reduce the computational costs of the method. Additionally, we conduct a study using different sets of snapshots, inspired by various approaches found in the literature. The snapshots have been computed from (i) the uncontrolled problem, (ii) the controlled problem, (iii) the linearized equation and also (iv) the  adjoint equation as suggested in \cite{SV13}.\\


The development of fast and accurate approximations for the SDRE helped to address our primary focus in the current work. Indeed, the integration of Proper Orthogonal Decomposition (POD) and Discrete Empirical Interpolation Method (DEIM) into the State-Dependent Riccati Equation (SDRE) framework, investigated in the first part of paper, is fundamental for the identification and control of our unknown system. At each time iteration, for the estimated parameter configuration, we solve a reduced Riccati equation and integrate it into the reduced black box, significantly expediting the procedure. In this context, given the high computational demand of solving the ARE, POD is of crucial importance. Finally, to update the estimation of the parameter configuration we employ Bayesian Linear Regression methods \cite{rossi2012bayesian}. Specifically, from the obtained trajectory, one can set a linear system that allows to estimate $\wt \mu$ at every time iteration. 
To the best of authors' knolwedge this is the first approach based on POD to identify and control a nonlinear system online.



It is noteworthy to acknowledge alternative approaches based on variants of sparse optimization techniques, such as Sparse Identification of Nonlinear Dynamics (SINDy, \cite{RBPK16,RABK19}), which has found application in both ODEs \cite{BPK16} and PDEs. SINDy has also been utilized for the identification of controlled systems, as demonstrated in \cite{KKB19}, where an external source was employed as input to identify the system, followed by the application of Nonlinear Model Predictive Control (NMPC) for system regulation. This ``identify first, control later" workflow differs from the strategy proposed in the current work. Additional strategies dedicated to the simultaneous control and identification of systems can be found in literature, such as the approach presented in \cite{KS08} for PDEs and \cite{MLG20} for ODEs. A recent study in \cite{CFDKMRV21} explores the control of unknown systems using Model Predictive Control (MPC), wherein the system identification is conducted using Extended Dynamic Mode Decomposition (EDMD), representing a surrogate linear model. This stands in contrast to our approach, where we directly identify the nonlinear model without relying on a surrogate linear representation.

To conclude we recap the novelties brought in the current contribution. First of all the use of POD in the identification of controlled problems using an online algorithm previously tested for low dimensional examples. Second, the use of POD-DEIM in the context of SDRE discussing different sets of snapshots. In both cases, the use of POD will show a computational benefit yielding a significant speedup factor along with the desired accurateness.

The outline of the paper is the following. In Section \ref{sec: sdre} we recall the State Dependent Riccati approach. In Section \ref{sec:podsdre}, we introduce the POD method for SDRE and provide a motivational example in Section \ref{sec:motivational}. In Section \ref{sec:idcon}, we describe all the details of the proposed method to identify and control an unknown problem using model reduction. Finally, in Section \ref{sec:test}
we show the effectiveness of our proposed method for two PDEs. Conclusions are driven in Section \ref{sec:conclusions}.

\section{Control of nonlinear problems via State Dependent Riccati Equations}\label{sec: sdre}
For the asymptotic stabilization of a nonlinear dynamics towards the origin, a suboptimal control can be obtained using the State Dependent Riccati Equation (SDRE) approach. We refer to e.g. \cite{C97,BLT07} for more details on this method.

Let $x(t):[0,\infty]\rightarrow\R^d$ be the state of the system, $u(t)\in\mathcal{U}:=L^{\infty}(\R^+;\R^m)$ be the control signal, $A(x):\R^d\rightarrow\R^{d\times d}$ and $B(x):\R^d\rightarrow\R^{d\times m}$ be state dependent matrices.

We consider a nonlinear dynamics in the form
\begin{align}\label{eq: dynamics}
\begin{aligned}
	\dot x(t)&=A(x(t))x(t)+B(x(t))u(t),\quad t\in(0,\infty),\\ x(0)&=x_0,
	\end{aligned}
\end{align}
associated to the following infinite horizon optimal control problem:
\begin{equation}\label{disc_cost}
\underset{u\in \mathcal{U}}{\min}\;J(u):=\int\limits_0^\infty \Big(\norm{x(t)}_Q^2+\norm{u(t)}_R^2\Big)\, dt
\end{equation}

where $Q\in\R^{d\times d},\,Q\succeq0$ and symmetric whereas $R\in\R^{m\times m},\,R\succ0$ and symmetric. Note that we will not deal with the most general case which considers the matrices $Q,R$ also state dependent.

Throughout this paper, the notation $\norm{v}_M^2$ stands for $v^{\top}Mv$ for any vector $v$ and compatible square matrix $M$. 

We can link, to this control problem, the state-dependent algebraic Riccati equation \begin{equation}\label{eq: sdre1}
A^\top(x)\Pi(x)+\Pi(x) A(x)-\Pi(x)B(x)R^{-1}B^\top(x)\Pi(x)+Q=0.
\end{equation}
We remark that, for a given state $x$, equation \eqref{eq: sdre1} is an algebraic Riccati equation (ARE). In what follows, we will refer to ARE also in the case of SDRE approach.
When this ARE admits a solution $\Pi(x)\in\R^{d\times d}$, we can define the feedback gain matrix as $K(x):=R^{-1}B^\top(x)\Pi(x)$. We refer to \cite{BLT07,BH18} for a detailed description on the assumption to guarantee existence of $\Pi(x)$. Then, a suboptimal nonlinear feedback law is given by
\begin{equation}\label{eq: sdref1}
u(x):=-K(x)x=-R^{-1}B^\top(x)\Pi(x)x.
\end{equation}

Thus, the idea is to solve an ARE at each time step in order to find a control. In \cite{BLT07} the following SDRE algorithm has been proposed.

\begin{algorithm}[ht]
	\caption{SDRE method }
        \label{alg: sdre}
 \begin{algorithmic}[1]
		\Require $\{t_0,t_1,\ldots\},$ model \eqref{eq: dynamics}, $R$ and $Q$,
		\For{$i=0,1,\ldots$}
		\State Compute $\Pi(x(t_i))$ from \eqref{eq: sdre1}
		\State Set $K(x(t_i)) := R^{-1}B^\top(x(t_i))\Pi(x(t_i))$
		\State Set $u(t):=-K(x(t_i))x(t), \quad \mbox{for }t \in [t_i, t_{i+1}]$ 
		\State Integrate the system dynamics with $u(t)$ to obtain $x(t_{i+1})$ \
		\EndFor
	\end{algorithmic}
\end{algorithm}

We remark that the ARE \eqref{eq: sdre1} does not always admit an analytical solution and the control found with this method is only suboptimal. In the current work, we approximate the solution of \eqref{eq: sdre1} using
rational Krylov subspaces (see \cite{S16}) if the dimension is large, e.g. $d\approx O(10^4)$, otherwise the Matlab function {\tt icare}.



\section{Proper Orthogonal Decomposition  for SDRE}\label{sec:podsdre}
The SDRE approach, discussed in Section \ref{sec: sdre}, becomes computationally expensive when dealing with a high-dimensional problem, such as $d \approx O(10^4)$. This implies solving the ARE in \eqref{eq: sdre1} at each iteration, involving $d^2$ unknowns.
 Therefore, it is useful to employ model order reduction techniques to reduce the computational cost of the control problem. 
 
 In this work, we specifically use the Proper Orthogonal Decomposition (POD, \cite{BHL93, S87}) method. While an initial approach to POD for SDRE can be found in \cite{HW23} for quadratic terms, we extend it to more general settings encompassing any nonlinear term. Additionally, a hyper-reduction method like Discrete Empirical Interpolation Method (DEIM, \cite{CS10,DG15}) will be employed.\\

We suppose to approximate $x(t)$ in \eqref{eq: dynamics} as $x(t)\approx \Psi x_r(t)$ where $\Psi\in\mathbb{R}^{d\times r}$ has orthonormal columns and $x_r(t):[0,\infty)\rightarrow\R^r$ with $r\ll d.$ If we plug our assumption into \eqref{eq: dynamics} and impose the Galerkin orthogonality we obtain the reduced dynamics
\begin{align}\label{eq: redynamics}
\begin{aligned}
	\dot x_r(t)&=A_r(x_r(t))x_r(t)+B_r(x_r(t))u(t),\quad t\in(0,\infty),\\ x_r(0)&=\Psi^T x_0,
	\end{aligned}
\end{align}
with $A_r(x_r):= \Psi^T A(\Psi x_r) \Psi\in\mathbb{R}^{r\times r},\, B_r(x_r) := \Psi^TB(\Psi x_r)\in\mathbb{R}^{r\times m}.$  The reduced cost functional we want to minimize reads
\begin{equation}\label{disc_cost}
J_r(u):=\int\limits_0^\infty \Big(\norm{x_r(t)}_{Q_r}^2+\norm{u(t)}_R^2\Big)\, dt
\end{equation}
with $Q_r = \Psi^T Q \Psi\in\mathbb{R}^{r\times r}$. 
Hence, for a given $x_r$, we obtain  the reduced ARE for the reduced problem \eqref{eq: redynamics}
\begin{equation}\label{eq: redsdre1}
A_r^\top(x_r)\Pi_r(x_r)+\Pi_r(x_r) A_r(x_r)-\Pi_r(x_r)B_r(x_r)R^{-1}B_r^\top(x_r)\Pi_r(x_r)+Q_r=0
\end{equation}
which is now a matrix equation for $\Pi_r(x_r)\in\R^{r\times r}$. The computational benefit of working with the small dimensional ARE \eqref{eq: redsdre1} instead of the one in \eqref{eq: sdre1} is clear. Finally, the reduced control will have the form
\begin{equation}\label{eq: sdref1red}
u_r(x_r):=-K_r(x_r)x_r=-R^{-1}B_r^\top(x_r)\Pi_r(x_r)x_r.
\end{equation}
Note that $u_r\in\mathbb{R}^m$ keeps the same dimension of the original problem. \\

We now discuss the selection of $\Psi$, a crucial aspect in POD applied to control problems. Typically, this method relies on snapshots (\cite{S87}) of the full-dimensional system. Let $x(t_i)$ be the (approximate) solution of \eqref{eq: dynamics} at time $t_i$ for a given input $u(t)$ and $\mathcal{S} = [x(t_0), x(t_1), \ldots, x(t_{n_t-1})]\in\R^{d\times n_t}$ be the snapshots matrix, i.e. a matrix of samples from the system \eqref{eq: dynamics}. The matrix $\Psi$, used to reduce the problem, is obtained from the left singular vectors of the rank $r$ Singular Value Decomposition (SVD) of $\mathcal{S}\approx \Psi \Sigma V^T.$  We will discuss various strategies to obtain the snapshots matrix $\mathcal{S}$, later in Section \ref{sec:motivational}. For control problems, it is well-known that the snapshots matrix should capture information relevant to the control problem. However, the optimal control input is typically unknown and snapshots obtained directly from the control problem are computationally expensive. For a comprehensive comparison of basis approaches in the context of feedback control, we refer to \cite{ASH17}.


\paragraph{Discrete Empirical Interpolation Method } We now turn our attention to the reduction of the nonlinear terms $A(x)$ and $B(x)$ since, in the actual reduced form, $A_r(x_r)$ and $B_r(x_r)$ still depend on the high dimensional problem through $\Psi x_r\in\R^d.$ To address this depency, we will employ the Discrete Empirical Interpolation Method (DEIM, \cite{CS10,DG15}) allowing the reduction independently of $d$. This is particularly effective when the nonlinear terms can be evaluated component-wise. Indeed, the DEIM method will select some components of the nonlinear terms to evaluate using a greedy strategy.

The DEIM method operates as follows: first, we compute the snapshots from the nonlinear term $\mathcal{F}_A = [A(x(t_0))x(t_0), A(x(t_1))x(t_1),\ldots A(x(t_{n_t-1}))x(t_{n_t-1})]\in\R^{d\times n_t}$. This set makes use of the snapshots set $\mathcal{S}$ previously computed. Next, we perform the reduced SVD of rank $\ell$ on $\mathcal{F}_A$, obtaining the left singular vectors denoted by $\Phi_A\in\mathbb{R}^{d\times\ell}.$ Those will be the POD basis for the nonlinear part. Then, we compute the QR factorization with pivoting of $\Phi_A^T$ and store the first $\ell$ columns of the permutation matrix into $P_A\in\mathbb{R}^{d\times\ell}$ obtained from the pivoting strategy.
Similarly for $B(x)$ one can follow the same approach with the snapshot set $\mathcal{F}_B = [B(x(t_0)), B(x(t_1)),\ldots B(x(t_{n_t-1})))]\in\R^{d\times (mn_t)}$. Therefore, one can use the DEIM approximation
\begin{align}
\begin{aligned}
A^\ell_r(x_r)&:=\Psi^T\Phi_A(P_A^T\Phi_A)^{-1} A(P_A^T\Psi x_r)P_A^T\Psi\approx A_r(x_r)\\
B^\ell_r(x_r)&:=\Psi^T\Phi_B(P_B^T\Phi_B)^{-1} B(P_B^T\Psi x_r)\approx B_r(x_r)
\end{aligned}
\end{align}
which are completely independent from the original dimension of the problem since the quantity $\Psi^T\Phi_A(P_A^T\Phi_A)^{-1}\in\mathbb{R}^{r \times \ell}$ and $P_A^T\Psi\in\mathbb{R}^\ell$ can be precomputed. Note that to distinguish the DEIM approach for $A(x)$ and $B(x)$ we use the index $A$ and $B$, respectively. One can also think to build a unique dataset $\mathcal{F}$ which includes information from all the nonlinear functions. This is not considered in the current work since our $B(x)$ will be constant in the numerical tests presented in Section \ref{sec:motivational} and later in Section \ref{sec:test}.
We will denote by $\Pi_r^\ell(x_r)$ the solution of the reduced ARE where it is used the $A_r^\ell(x_r)$ instead of $A_r(x_r)$ in \eqref{eq: redsdre1} which reads
\begin{equation}\label{eq: poddieimsdre}
(A^\ell_r)^\top(x_r)\Pi_r^\ell(x_r)+\Pi_r^\ell(x_r) A^\ell_r(x_r)-\Pi_r^\ell(x_r)B^\ell_r(x_r)R^{-1}(B^\ell_r)^\top(x_r)\Pi_r^\ell(x_r)+Q_r=0
\end{equation}
For completeness we also provide the POD-DEIM dynamics:
\begin{align}\label{eq: deimdyn}
\begin{aligned}
	\dot x_r(t)&=A_r^\ell(x_r(t))x_r(t)+B^\ell_r(x_r(t))u(t),\quad t\in(0,\infty),\\ x_r(0)&=\Psi^T x_0.
	\end{aligned}
\end{align}

The DEIM approach is well-known in the framework of model reduction and we refer to \cite{DG15} for extensive details. However, to the best of the authors' knowledge it is applied for the fist time to the SDRE problem. To summarize we show the POD-DEIM method for SDRE in Algorithm \ref{alg: podsdre}.

\begin{algorithm}[ht]
	\caption{POD-DEIM for SDRE method }
        \label{alg: podsdre}
 \begin{algorithmic}[1]
		\Require $\{t_0,t_1,\ldots\},$ model \eqref{eq: dynamics}, $R$ and $Q$, POD base $\Psi$, DEIM base $\Phi$, matrix $P$
		\For{$i=0,1,\ldots$}
		\State Compute $\Pi_r^\ell(x(t_i))$ from \eqref{eq: poddieimsdre}  
		\State Set $K_r(x(t_i)) := R^{-1}(B^\ell_r)^\top(x_r(t_i))\Pi_r^\ell(x(t_i))$
		\State Set $u_r(t):=-K_r(x_r(t_i))x_r(t), \quad \mbox{for }t \in [t_i, t_{i+1}]$ 
		\State Integrate the reduced system dynamics \eqref{eq: deimdyn} with $u_r(t)$ to obtain $x_r(t_{i+1})$ 
		\EndFor
   \State Obtain $x(t)\approx \Psi x_r(t)$

	\end{algorithmic}
\end{algorithm}
We remark that the approximate solution of the reduced ARE in \eqref{eq: poddieimsdre} is computed using the Matlab function {\tt icare} due to the low dimensionality.

\subsection{Numerical Example: control of the 2D Allen-Cahn}\label{sec:motivational}

In this section we compare the results of Algorithm \ref{alg: sdre} and Algorithm \ref{alg: podsdre} for the control of the 2D Allen-Cahn equation. Specifically, we will study the accurateness and the computational costs of the reduced approach for different snapshot sets.

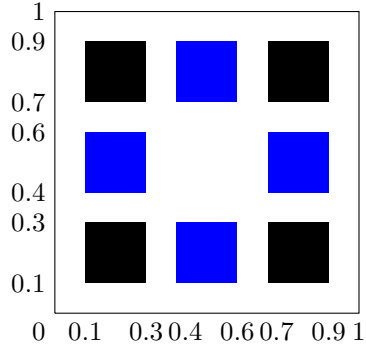
\begin{figure}[htbp]
\centering
\begin{tikzpicture}[scale=4]
\draw (0,0) -- (1,0) -- (1,1) -- (0,1) -- (0,0); 
\fill[black] (0.1,0.1) -- (0.3,0.1) -- (0.3,0.3) -- (0.1,0.3) -- (0.1,0.1);
\fill[black] (0.7,0.7) -- (0.9,0.7) -- (0.9,0.9) -- (0.7,0.9) -- (0.7,0.7);
\fill[black] (0.1,0.7) -- (0.3,0.7) -- (0.3,0.9) -- (0.1,0.9) -- (0.1,0.7);
\fill[black] (0.7,0.1) -- (0.9,0.1) -- (0.9,0.3) -- (0.7,0.3) -- (0.7,0.1);
\fill[blue]  (0.1,0.4) -- (0.3,0.4) -- (0.3,0.6) -- (0.1,0.6) -- (0.1,0.4);
\fill[blue] (0.1,0.4) -- (0.3,0.4) -- (0.3,0.6) -- (0.1,0.6) -- (0.1,0.4);
\fill[blue] (0.7,0.4) -- (0.9,0.4) -- (0.9,0.6) -- (0.7,0.6) -- (0.7,0.4);
\fill[blue] (0.4,0.1) -- (0.6,0.1) -- (0.6,0.3) -- (0.4,0.3) -- (0.4,0.1);
\fill[blue] (0.4,0.7) -- (0.6,0.7) -- (0.6,0.9) -- (0.4,0.9) -- (0.4,0.7);
\draw (0,0) node[below left] {$0$};
\draw (0.1,0) node[below ] {$0.1$};
\draw (0.3,0) node[below ] {$0.3$};
\draw (0.4,0) node[below ] {$\,\,\,\,0.4$};
\draw (0.6,0) node[below ] {$0.6$};
\draw (0.7,0) node[below ] {$\,\,\,\,0.7$};
\draw (0.9,0) node[below ] {$0.9$};
\draw (1,0) node[below ] {$1$};
\draw (0,0.1) node[left ] {$0.1$};
\draw (0,0.3) node[left ] {$0.3$};
\draw (0,0.4) node[left ] {$0.4$};
\draw (0,0.6) node[left] {$0.6$};
\draw (0,0.7) node[left] {$0.7$};
\draw (0,0.9) node[left ] {$0.9$};
\draw (0,1) node[left ] {$1$};
\end{tikzpicture}
\caption{The control, or input, acts in the black region $\Omega_c$ and the output, i.e. the cost, is considered in the blue region $\Omega_o=\bigcup_{i=1}^z\Omega_{o_i}$.}
\label{fig:bc_coll}
\end{figure}

We consider the control of a two-dimensional semilinear parabolic equation over $\Omega\times\R^+_0$, with $\Omega=(0,1)^2\subset \R^2$ and Dirichlet boundary conditions:

\begin{equation}\label{ex1}
    \left\{ \begin{aligned}
    y_t(t,\xi) &= \mu_1 \Delta y(t,\xi) + \mu_2y(t,\xi) + \mu_3 y^3(t,\xi) + \chi_{\Omega_c}(\xi)u(t),\\
    y(t,\xi) &= 0, \quad \xi\in\partial\Omega,\\
    y(0,\xi) &=y_0(\xi), \quad \xi\in\Omega
    \end{aligned} \right.
\end{equation}
with $\mu_1=0.5,\mu_2=11,\mu_3=-11$.
The scalar control $u(t)$ acts through the indicator function $\chi_{\Omega_c}(\xi)$  with support $\Omega_c\subset\Omega$ defined as
$$
\Omega_c = [0.1,0.3]^2 \cup [0.7,0.9]^2 \cup 
\left([0.1,0.3]\times[0.7,0.9]\right) \cup \left( [0.7,0.9]\times[0.1,0.3]\right)
$$
The cost, taken from \cite{AKS23}, is given by
\begin{align}\label{costex1}
\begin{aligned}
J(u)&:=\int\limits_0^{\infty} \sum_{i=1}^z \dfrac{1}{|\Omega_{o_i}|}\left(\int_{\Omega_{o_i}} y(t,\xi)\, d\xi\right)^2 +R|u(t)|^2  \, dt\,
\end{aligned}
\end{align}
where $z=4$ and
\begin{align*}
\begin{aligned}
\Omega_{o_1}  =& [0.1,0.3]\times[0.4,0.6],\quad \Omega_{o_2}  =& [0.4,0.6]\times[0.1,0.3]\\
\Omega_{o_3}  =& [0.4,0.6]\times[0.7,0.9],\quad \Omega_{o_4}  =& [0.7,0.9]\times[0.4,0.6].
\end{aligned}
\end{align*}
Therefore, the control and the cost act in different regions of $\Omega$, as shown in Fig.  \ref{fig:bc_coll}.

Equation \eqref{ex1} admits a trivial solution $y(t)\equiv 0$ which is an unstable equilibrium.
We notice that for the chosen model and $\mu$ parameter's values, the uncontrolled solution ($u\equiv 0$) does not tend to that equilibrium
point. Indeed, the aim of our control problem is to stabilize the system towards the equilibrium $y(t)\equiv 0$. 

To set equation \eqref{ex1} into the form \eqref{eq: dynamics}, we discretize in space the system dynamics, using e.g. a finite difference approximation, and write it in a semi discrete form.
We define the discrete state $x(t)=(x_1(t),\ldots,x_d(t))^\top\in\R^d$ as the approximation of $y(t,\xi)$ at the grid points. Therefore, if there are respectively $n_{\xi_1}$ and $n_{\xi_2}$ points in the discretization of $\xi_1$ and $\xi_2$ axis, then $d=n_{\xi_1}n_{\xi_2}$.


%
The semi discrete equation from \eqref{ex1} reads
\begin{equation}\label{pol:sd}
\dot x(t) = \Bigl(\mu_1\Delta_d+\mu_2\textbf{1}_{d}+\mu_3 \text{diag}(x^2(t))\Bigr)x(t) +B u(t)
\end{equation}
where $\Delta_d\in\R^{d\times d}$ denotes the second order finite difference approximation of the  
Dirichlet Laplacian, $\text{diag}(x(t)^2)\in \R^{d\times d}$ indicates a diagonal matrix with the squares of the elements of $x(t)$ on the diagonal, $\textbf{1}_{d}\in \R^{d\times d}$ is the $d$-dimensional identity matrix,  and $B\in \R^d$ is the discretization of the indicator function supported over $\Omega_c$. Hence, we obtain a problem in the form \eqref{eq: dynamics} where \mbox{$A(x)=\mu_1\Delta_d+\mu_2\textbf{1}_{d}+\mu_3 \text{diag}(x^2)$}.
The discretization of the cost \eqref{costex1} can be written in the form \eqref{disc_cost} by properly setting $Q$ (see \cite{AKS23}).

We consider the initial condition $y(0,\xi) = 0.2\sin (\pi\xi_1)\sin(\pi\xi_2)$, on a discretized space grid of $n_{\xi_1} n_{\xi_2}$ nodes with $n_{\xi_1} = n_{\xi_2} = 101$ ($d=10201$). The time discretization will be performed with a step size $\Delta t =0.025$ and integrated using an implicit Euler scheme. We also note that, although the problem deals with an infinite horizon, in the numerical simulations we chose an horizon large enough, say $T=3$ in this example. The reported numerical simulations were performed on a iMac with Apple M1 and 16GB RAM, using Matlab \cite{matlab7}.

\begin{figure}[htbp]
\centering
\includegraphics[scale=0.3]{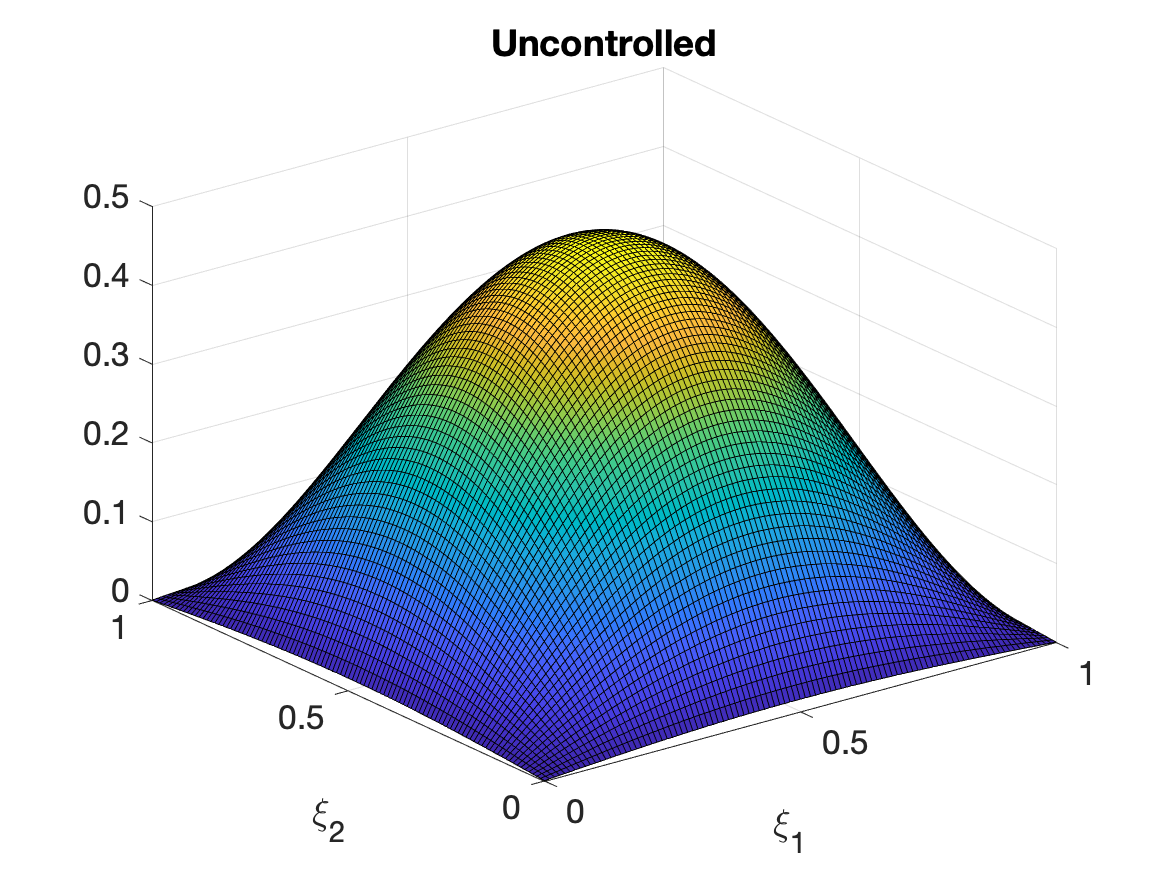}
\includegraphics[scale=0.3]{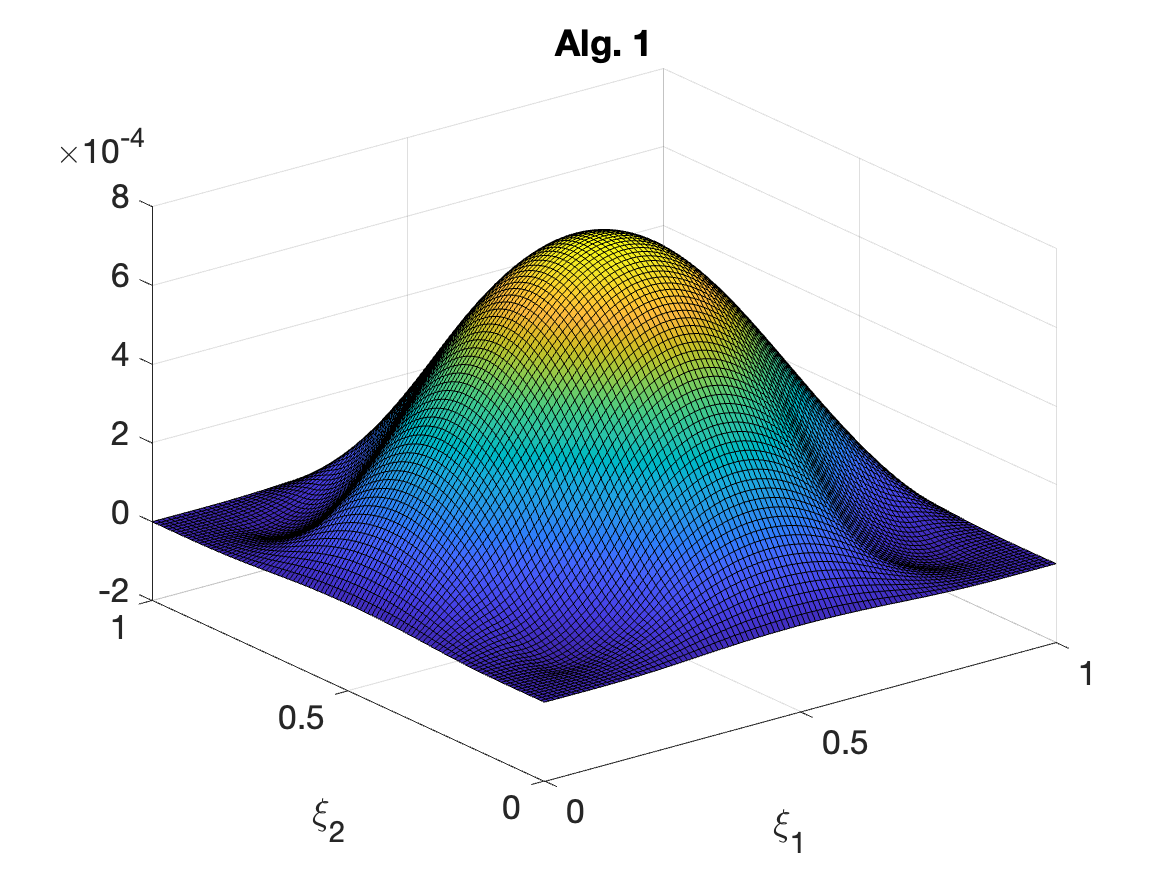}
\caption{Solutions at time $t=3$ of \eqref{ex1} for the uncontrolled problem $u(t)=0$ (left), and the stabilized solution through Algorithm \ref{alg: sdre} (right). }
\label{fig:ex1}
\end{figure}

In Figure \ref{fig:ex1}, we compare the approximate solutions at time $t=3$ for the uncontrolled problem, i.e. $u(t)\equiv 0$ in \eqref{ex1} in the left panel, whereas we show the stabilized solution using Algorithm \ref{alg: sdre} for \eqref{ex1}. It is visually clear, from the scaling of the $z-$axis, the difference and the effectiveness of the algorithm to stabilize the solution. The CPU time to perform Algorithm \ref{alg: sdre} was about 135s.

Next, we want to show the performance of the POD-DEIM method for SDRE proposed in  Algorithm \ref{alg: podsdre}.
We will provide 4 different approaches used to compute the snapshots matrix. We will denote by $x(t_i;u,\mu)$
the approximate trajectory at time $t_i$ using control $u$ with the parameter configuration $\mu$. This will be useful to stress which control or parameter configuration will be used to compute the snapshots. If we use the configuration $\mu_1 = 0.2, \mu_2 = 11, \mu_3 = -11,$ we will drop the depency on $\mu$ and use the notation $x(t_i;u)$. Below, we provide the choices of snapshots we investigated in this section:

\begin{enumerate}
    \item First approach consists in choosing the snapshots from the uncontrolled problem, i.e. using $u\equiv 0$. This will be denoted by  $\mathcal{S}_0=[x(t_0;0),\ldots, x(t_{n_t-1};0)].$
    \item Second approach consists in choosing the snapshots from the solution of Algorithm \ref{alg: sdre}: $\mathcal{S}_{u^*}=[x(t_0; u^*(t_0)),\ldots, x(t_{n_t-1};u^*(t_{n_t-1}))].$
    \item Third approach consists in computing the optimal control of the linearized problem. The linearization of \eqref{ex1} around the equilibrium $y(t,\xi) = 0$, which is the goal in our optimal control problem, corresponds to the same model with $\mu_3=0$. Its semi discrete form is an LQR problem (see \cite{anderson1999LQR}):
   
    \begin{equation}\label{pol:sd_lin}
        \dot x_{lin}(t) = \Bigl(\mu_1\Delta_d+\mu_2\textbf{1}_{d}\Bigr)x_{lin}(t) +B u_{lin}(t).
    \end{equation}
 Furthermore, since in the sequel our problem will be only partially known, we will assume that the parameters in \eqref{pol:sd_lin} belong to certain intervals: $\mu_1\in I_1\subset\R$ and $\mu_2\in I_2\subset\R.$ We, then, discretize the intervals $I_1, I_2$, thus obtaining a finite number $n_{par}$ of possible values $\{\hat \mu_1,\ldots, \hat\mu_{n_{par}}\}.$ Note that in this case $\hat\mu_j\in\R^2.$
    For each combination $j\in\{1,\ldots ,n_{par}\}$, we compute the LQR control $u_{lin}(t;\hat \mu_j)$\footnote{The notation $u_{lin}(t;\hat \mu_j)$ stands for the control from \eqref{pol:sd_lin} with parameter configuration $\hat\mu_j$. In practice, we compute the feedback gain matrix $K_j^{lin}$ and set $u^{lin}(t;\hat\mu_j)=-K_j^{lin}x(t)$ where $x(t)$ relates to \eqref{pol:sd}.} and form a snapshots set from \eqref{pol:sd}. We set  $\mathcal{S}_{j}=[x(t_0; u_{lin}(t;\hat\mu_j)); \ldots; x(t_{n_t-1}; u_{lin}(t;\hat\mu_j))]$. The considered snapshots set will be $\mathcal{S}_{lin}:=\cup_{j=1}^{n_{par}} \mathcal{S}_{j}$. In other words, here, the snapshots are computed from the original problem \eqref{pol:sd} using the feedback gain obtained from its linearized equation \eqref{pol:sd_lin} for different parameter configurations. The advantage of the linearized approach is the independence of the associated Riccati equation from the state. This allows to solve only one ARE and clearly reduces the cost of the snapshots computation.
    \item Similarly to the previous step we also provide adjoint information to the snapshots set motivated by the work in \cite{SV13}. For completeness, we provide the adjoint problem related to the linearized version of \eqref{ex1}
    \begin{equation}\label{adj}
    \left\{ \begin{aligned}
    -p_t(t,\xi) &= \mu_1 \Delta p(t,\xi) + \mu_2 p(t,\xi) -y(t,\xi)\\
    p(t,\xi) &= 0, \quad \xi\in\partial\Omega\\
    p(T,\xi) &=0, \quad \xi\in\Omega.
    \end{aligned} \right.
\end{equation}
Its semi discrete form for $\tilde p(t)\approx p(t,\xi)$ reads:
\begin{equation}\label{pol:adj_lin}
       - \dot {\tilde p}(t) = \Bigl(\mu_1\Delta_d+\mu_2\textbf{1}_{d}\Bigr)\tilde p(t) -x_{lin}(t).
    \end{equation}
The snapshots from the adjoint problem are denoted by $\mathcal{S}_{adj} := \cup_{j=1}^{n_{par}} \mathcal{P}_j$ with $\mathcal{P}_j=[\tilde p(t_0;\hat\mu_j),\ldots, \tilde p(t_{n_t};\hat\mu_j)]$.  We will then investigate, as forth approach, $\mathcal{S}_{lin}\cup\mathcal{S}_{adj}$. The notation $p(t_i;\hat\mu_j)$ stresses the dependence on the parameter configuration $\hat\mu_j$.\\

\end{enumerate}

\begin{remark}
The relevance of the third approach will be clearer, later in Section \ref{sec:idcon}. There, we study a model which is not completely know, and therefore it will be important to set a reduced model which surrogates many parametric configurations.\\
\end{remark}



To measure the quality of the POD-DEIM approximation, we will use the relative error $\mathcal{E}(r)$ between the SDRE solution from Algorithm \ref{alg: sdre}, which is used as ``reference" (or exact) solution, and its POD approximation from Algorithm \ref{alg: podsdre} and the difference $\mathcal{E}_J(r)$ between the cost functionals:
$$\mathcal{E}(r)=\max_t\frac{\|x(t) - \Psi x_r(t) \|_2}{\|x(t)\|_2},\qquad \mathcal{E}_J(r)=|J(u)-J_r(u_r)|. $$

\begin{figure}[hbtp]
\includegraphics[scale=0.25]{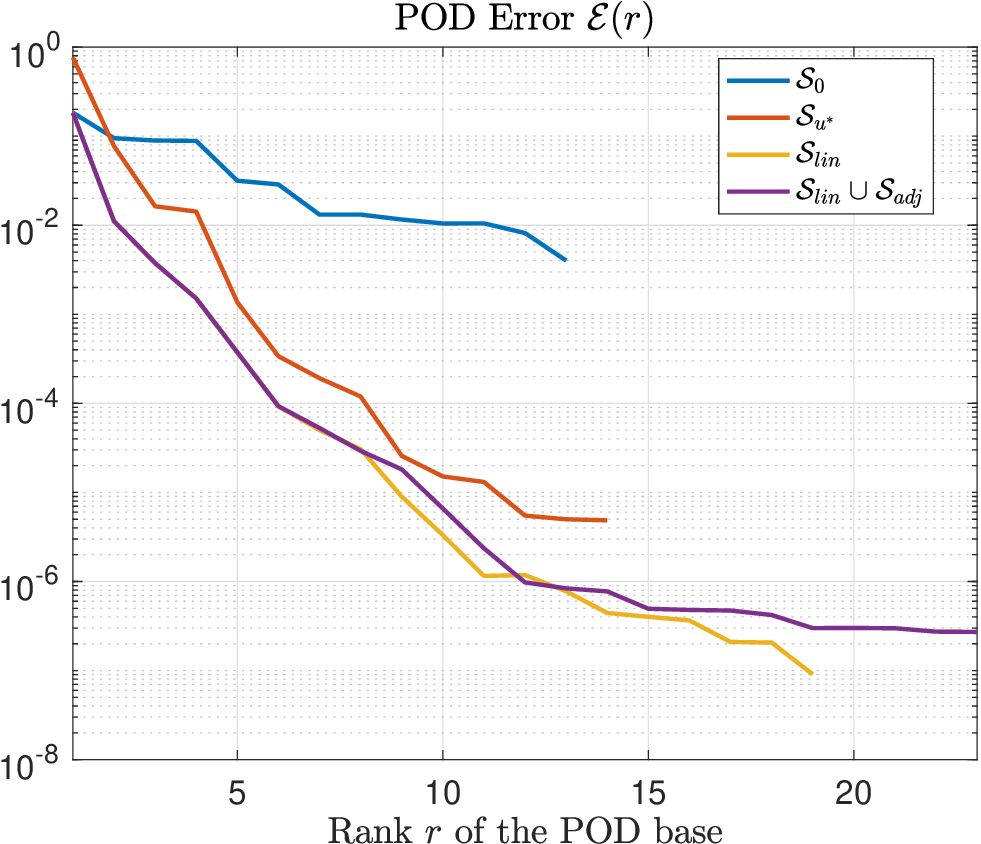}
\includegraphics[scale=0.25]{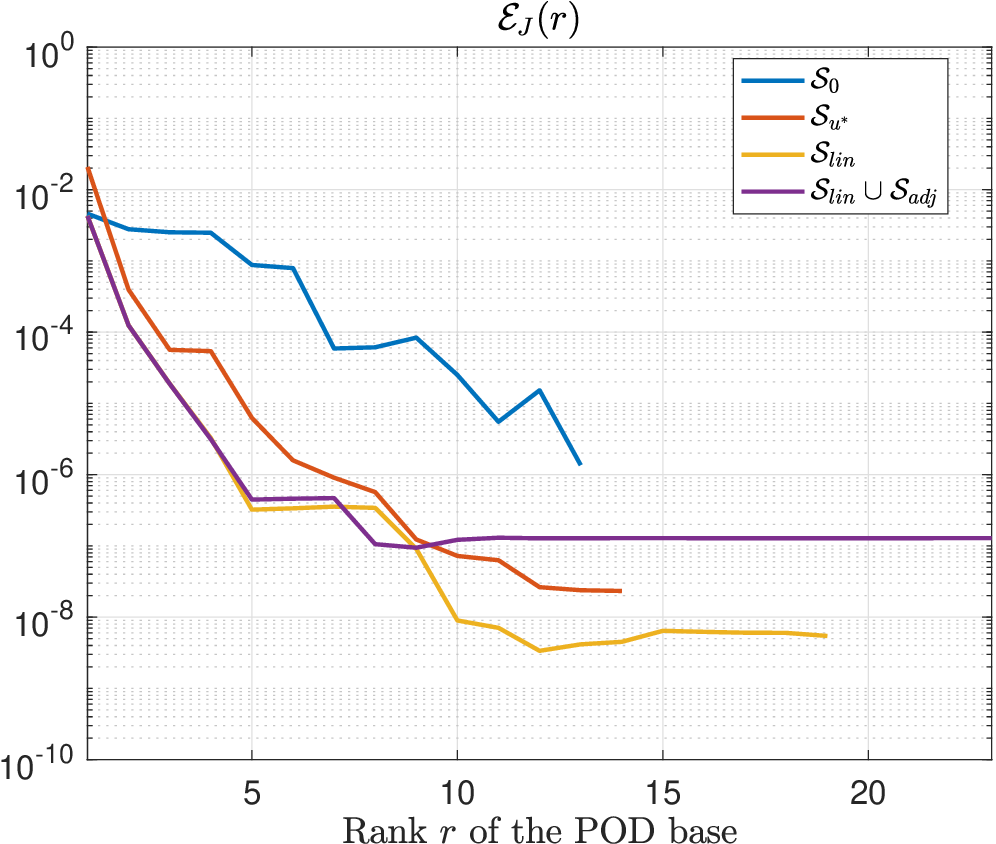}
\includegraphics[scale=0.25]{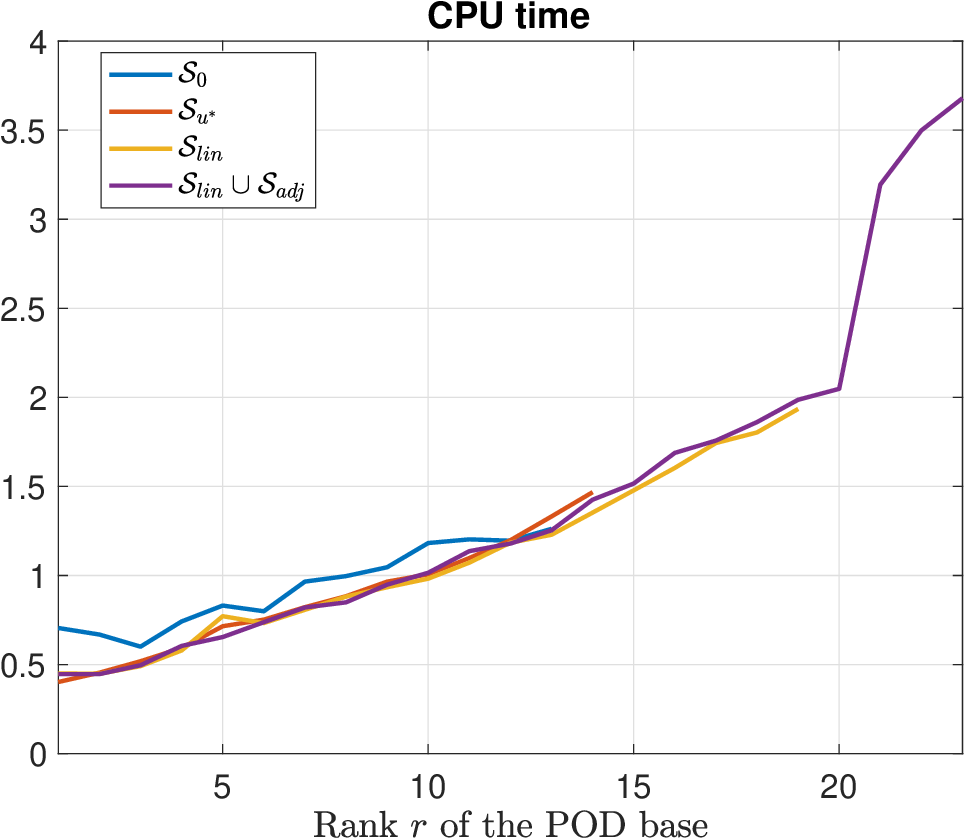}
\caption{Results of Algorithm \ref{alg: podsdre} for different snapshots set. POD error $\mathcal{E}(r)$ (left), Difference for the cost functional $\mathcal{E}_j(r)$ (middle), CPU time (right). The time needed to solve Alg. \ref{alg: sdre} was 135s. We can observe an impressive speed up of the POD-DEIM method.}
\label{fig:pod1}
\end{figure}

In Figure \ref{fig:pod1}, we show the results of POD-DEIM applied to the control problem of \eqref{ex1}. For the third and fourth approach we considered a discretization of the intervals $I_1, I_2$ considering $\{0.1, 0.2154, 0.4642 \}\subset I_1 $ and $\{0, 7.5, 15\}\subset I_2$. Note that the parameter configuration studied for \eqref{ex1} is not included in the snapshots set. The number of considered snapshots $n_t$ is $n_t = 120$ for $\mathcal{S}_0$ and $\mathcal{S}_{u^*}$, $n_t=270$ for $\mathcal{S}_{lin}$ and $n_t = 540$ for $\mathcal{S}_{lin}\cup\mathcal{S}_{adj}.$

The behavior of $\mathcal{E}(r)$ is shown in the left panel. We can see the decay of the error when the number of POD basis $r$ is increased. The number of DEIM points $\ell$ is fixed to the rank of $\mathcal{F}_A$ which in our simulations was $\ell=13$ for the snapshots related to $\mathcal{S}_0$ and $\mathcal{S}_{u^*}$, $\ell=20$ for $\mathcal{S}_{lin}$ and $\ell=21$ for $\mathcal{S}_{lin}\cup\mathcal{S}_{adj}$. One can see in the left panel of Figure \ref{fig:pod1} that the set of snapshots $\mathcal{S}_0$ does not provide very accurate solutions, while the sets $\mathcal{S}_{lin}$ and $\mathcal{S}_{lin}\cup\mathcal{S}_{adj}$ produce very similar accurate performances. For this reason, in the following, we will use $\mathcal{S}_{lin}$ to select our snapshots set. This choice is also confirmed by the absolute cost difference in $\mathcal{E}_j(r)$ (see middle panel of Fig \ref{fig:pod1}). Finally, we show the computational benefit of the POD strategies in the right panel of Fig \ref{fig:pod1}. It is also clear that the augmented set of snapshots which includes adjoint information requires a larger CPU time. This further motivates our choice for the snapshots set in what follows.

\section{Online Identification and Control of Unknown Nonlinear Dynamics using POD}\label{sec:idcon}

Our objective is to regulate a partially unknown nonlinear system by unraveling its characteristics online through continuous trajectory observation. This is achieved by providing an initial condition together with an input
$u(t)$ to obtain the resulting trajectory. The system under consideration, derived from \eqref{eq: dynamics}, is described by:
\begin{align}\label{eq:dynmu}
\begin{aligned}
	\dot x(t) &= \sum_{j=1}^n \mu_j A_j(x(t))x(t) + B(x(t))u(t),\quad t\in(0,\infty),\\
	x(0) &= x_0,
\end{aligned}
\end{align}
with $A(x) = \sum_{j=1}^n \mu_j A_j(x)$ and $A_j(x):\R^d\rightarrow\R^{d\times d}$ for $j= 1,\ldots,n$. The coefficients $\mu_j$ allow flexibility in selecting terms from the known library $A_j(x)$. 
Our problem consists in estimating the parameter configuration $\mu=(\mu_1,\ldots,\mu_n)\in\mathcal{D}\subset\R^n$ with $\mathcal{D}$ beging a compact set and computing a control $u(t)$ that stabilizes the unknown model. The unknown model relies on a parameter configuration which will be denoted by $\mu^*=(\mu_1^*,\ldots,\mu_n^*)\in\R^n$ and it is observable only through the trajectories in \eqref{eq:dynmu}. Furthermore, in this work, we assume that we can observe the whole system to compute the snapshots. Then, the system will be observed only once during the online control and identification phase.

\noindent
Let us, first, fix the notations useful throughout this section:

\begin{itemize}
\item $\mu$ is a generic parameter,
\item$\mu^*$ is the true system configuration,
\item $\wt\mu$ is an estimate system configuration,
\item $\wt\mu_r$ is an estimate system configuration obtained using model reduction,
\item $u(t;\wt \mu)$ is the control for the estimated configuration $\wt\mu$ computed through e.g. Algorithm \ref{alg: sdre},
\item $x(t; u(t;\wt \mu),\mu^*)$ is the observed trajectory computed with the input $u(t;\wt \mu)$ related to the true model $\mu^*,$
\item $x^i(u(t_i;\wt\mu),\mu^*) = x(t_i;u(t_i;\wt \mu),\mu^*)$ is the above-mentioned solution at discrete time $t_i$ where the control $u(t_i;\wt \mu)$ is applied in the interval $[t_i,t_{i+1}]$
\item $u_r(t;\wt \mu_r)$ is the control for the estimated configuration $\wt\mu_r$ computed through e.g. Algorithm \ref{alg: podsdre},
\item $x_r(t; u_r(t;\wt \mu_r),\mu^*)$ is the observed reduced trajectory computed with the input $u_r(t;\wt \mu_r)$ related to the true model $\mu^*,$
\item $x_r^i(u_r(t_i;\wt\mu_r),\mu^*) = x_r(t_i;u_r(t_i;\wt \mu_r),\mu^*)$ is the above-mentioned solution at discrete time $t_i$ where the control $u_r(t_i;\wt \mu_r)$ is applied in the interval $[t_i,t_{i+1}]$.
 \end{itemize}




The cost functional to minimize is adapted from \eqref{disc_cost} and reads:
\begin{equation}\label{disc_costmu}
J(u;\mu^*) = \int\limits_0^\infty \Big(\norm{x(t; u,\mu^*)}_Q^2+\norm{u(t)}_R^2\Big)\, dt,
\end{equation}
where $x$ is observed from the true system configuration for a given input.

In the context of the system \eqref{eq:dynmu}, the focus is on two unknowns: (i) the parameter configuration $\mu$, which must converge to an unknown true configuration $\mu^*$, and (ii) the control $u(t)$, that should minimize \eqref{disc_costmu} and stabilize the trajectory in \eqref{eq:dynmu}. This addresses the joint problem of system identification and control for \eqref{eq:dynmu}. The computation of the control is executed using a SDRE approach, according to a particular parameter configuration, that is approximated at each iteration using system observations. In \cite{RLPDE}, the authors have presented the results for small scale problems and shown their numerical convergence. The algorithm will be recalled later in this section in Algorithm \ref{alg:RL}.\\

Throughout this paper we want to  tackle more challenging problems since we will deal with two-dimensional PDEs. This requires a keen focus on the computational efficiency of the method. As already shown in Section \ref{sec: sdre}, the computational cost of Algorithm \ref{alg: sdre} is very high, and for this reason we have introduced model reduction to enhance efficiency. In what follows we will provide a fully reduced algorithm to achieve our goal. 



Our proposed method follows a classical offline/online decomposition which is typical in the model reduction framework. The offline stage includes the computation of the snapshots and of all the quantities to set the reduced problem. This task is usually computationally expensive, but it is performed only once at the beginning. The online stage, instead, includes at each time step the computation of the control, the trajectory observation and the parameter update which will be explained in detail. The computation of the control is performed using the reduced system, and this allows to avoid solving SDRE in high dimension using Algorithm \ref{alg: podsdre}. The trajectory observation will be carried out from the reduced system. The parameter update is carried out using BLR in the full system since it is well-known that the POD dynamics does not preserve the original structure of the problem. The algorithm is described in Algorithm \ref{alg:RL POD red}.

\begin{algorithm}[H]
	\caption{POD identification and control}
	\label{alg:RL POD red}
	\begin{algorithmic}[1]
		\Require $\{t_0,t_1,\ldots\},$ model $\{A_j(x)\}_{j=1}^n, B, R,Q$, $\wt{\mu}^0_r$, $\Sigma_0$
        \State Compute the snapshots \label{step: snap}
        \State Project the problem \label{step: proj}
		\For{$i=0,1,\ldots$}
			\State Solve \eqref{eq: redsdre1} and obtain $\Pi_r^{\ell}(x_r(t_i);\wt\mu_r^i)$\label{step: pi}
			\State Set $K_r(x_r(t_i);\wt{\mu}^i_r)$ using \eqref{eq: sdref1red}\label{step: K}
			\State Set $u_r(t;\wt{\mu}^i_r):=-K_r(x_r(t_i);\wt{\mu}^i_r)x_r(t)$\quad $\mbox{for }t \in [t_i, t_{i+1}]$ \label{step: u}
			\State Observe the trajectories $x_r^{i+1}(u_r(t_i;\wt \mu^i_r),\mu^*)$\label{step: traj}
                \State Compute $X_r$ and $Y_r$ as in \eqref{eq: XY BLR}\label{step: BLR1}
                \If{$\|\wt{\mu}^{i}_r-\wt{\mu}^{i-1}_r\|_\infty<tol_\mu$}
			      \State Compute $\wt{\mu}^{i+1}_r$ using BLR from \eqref{eq: BLR formulas}\label{step: BLR2}
                \Else
                    \State $\wt{\mu}^{i+1}_r=\wt{\mu}^{i}_r$\label{step: BLRcond}
                \EndIf               
		\EndFor
	\end{algorithmic}
\end{algorithm}

Let us now comment and describe each step of Algorithm \ref{alg:RL POD red}.\\

\noindent
{\bf Offline stage:}

\paragraph{ Choose an initial parameter configuration.}  We provide an initial  probability distribution for the parameter $\mu\sim\mathcal{N}(\wt{\mu}^0_r,\Sigma_0)$, that will be used as prior distribution
 for the true system configuration $\mu^*$. The initial parameter estimate at time $t=0$ is $\wt{\mu}^0_r\in\R^n$ whereas $\Sigma_0\in\mathbb{R}^{n\times n}$ is the covariance matrix.  If prior information about $\mu^*$ is available, it can guide the choice of $\wt{\mu}^0_r$. In the numerical tests, we will set $(\wt{\mu}^0_r)_k=1$ for $k=1,\ldots n$, and  use a heuristic covariance matrix $\Sigma_0 = c \textbf{1}_n$ with $c>0$ large. Later, the notation $\wt\mu^i_r$ will denote the parameter estimate, i.e. the mean of the obtained posterior distribution, at time $t_i$.
    \paragraph{ Compute the snapshots (step \ref{step: snap}). } Our choice is the use of $\mathcal{S}_{lin}$ discussed in Section \ref{sec:motivational}.
   \paragraph{ Project the problem into the reduced subspace (step \ref{step: proj}).}  Here we compute all the quantities needed to set our reduced problem as for instance $\Psi, \Phi, P$. If there are some linear terms such that for some index $\bar j$, the matrix $A_{\bar j}$ is constant we can also set its POD projection $\Psi^T A_{\bar j} \Psi.$\\

\noindent
{\bf Online stage:} we repeat at each time step $t_i$:
\paragraph{ Compute the reduced control (steps \ref{step: pi}-\ref{step: u}).} 
    At time $t_i$, we compute an approximate solution for the reduced ARE \eqref{eq: redsdre1} obtained from the current parameter estimate $\wt{\mu}^i_r$, thus obtaining the reduced feedback gain $K_r$ and the reduced feedback control $u_r(t;\wt{\mu}^i_r)$. For the first time step ($i=0$) we set $u=0$ because the first parameter approximation $\wt\mu^0$ can be far from the true $\mu^*$.

\paragraph{ Observe the trajectories (step \ref{step: traj}).}
    We apply piecewise-constant control $u_r(t;\wt{\mu}^i_r)$ in the time interval $[t_i,t_{i+1}]$ and observe the reduced trajectory $x_r^{i+1}(u_r(t_i;\wt\mu^i_r);\mu^*)$.
    
   
 \paragraph{ Update the parameter configuration based on the observations (steps \ref{step: BLR1}-\ref{step: BLRcond}).}
 We update the parameter estimate $\wt\mu^{i+1}_r$ using BLR (see e.g. \cite{Rasmussen2006}). In order to do this, we first assume $\mu^*\sim \mathcal{N}(\wt\mu^i_r,\Sigma_i)$, where $\wt\mu^i_r$ and $\Sigma_i$ are known from the previous time step. Discretizing the system \eqref{eq:dynmu}, we obtain a linear system of equations in the form $Y_r=X_r\mu^*+\varepsilon$, where $\mu^*$ is the true but unknown parameter configuration, $\varepsilon$ is a noise representing the error introduced with the discretization, that is assumed to be gaussian $\varepsilon\sim \mathcal{N}(0,\sigma^2)$, and $X_r\in \mathbb{R}^{d\times n}$ and $Y_r\in \mathbb{R}^d$ are given by 
\begin{align}
\begin{aligned}
    \label{eq: XY BLR}
        X_r&= [A_1(\Psi x_r^{i+1})\Psi x_r^{i+1},\ldots, A_n(\Psi x_r^{i+1})\Psi x_r^{i+1}],\\
        Y_r&=\frac{\Psi x_r^{i+1}-\Psi x_r^i}{\D t} - B u_r(t_i;\wt\mu^i_r) 
    \end{aligned}
    \end{align}  
    where $x_r^{i+1}=x_r^{i+1}(u_r(t_i;\wt\mu^i_r);\mu^*)$. 
    Therefore, using BLR we obtain a posterior distribution $\mu^*\sim\mathcal{N}(\wt\mu^{i+1}_r,\Sigma_{i+1})$, where $\wt\mu^{i+1}_r$ and $\Sigma_{i+1}$ can be computed explicitly as
    \begin{equation}\label{eq: BLR formulas}
        \Sigma_{i+1} =\biggl( \frac{1}{\sigma^2} X_r^T X_r + \Sigma_i^{-1}\biggr)^{-1}, \quad \wt\mu^{i+1}_r = \Sigma_{i+1} \left( \frac{1}{\sigma^2} X_r^T Y_r + \Sigma_i^{-1} \wt\mu^i_r \right).
    \end{equation}\\

For completeness we recall in Algorithm \ref{alg:RL} the method proposed in \cite{RLPDE} where model reduction was not involved and that was only applied to one dimensional PDEs. We will use it in the numerical tests to compare the results obtained with Algorithm \ref{alg:RL POD red}.

The expressions of $X$ and $Y$ in this case are
\begin{align}
\begin{aligned}
    \label{eq: XY_full}
        X&= [A_1(x^{i+1})x^{i+1},\ldots, A_n(x^{i+1})x^{i+1}],\\
        Y&=\frac{ x^{i+1}-x^i}{\D t} - B u(t_i;\wt\mu^i) 
    \end{aligned}
    \end{align}
while $\wt\mu^{i+1}$ is updated using formulas similar to \eqref{eq: BLR formulas}.

\begin{algorithm}[H]
	\caption{Online Identification and Control from \cite{RLPDE}}
	\label{alg:RL}
	\begin{algorithmic}[1]
		\Require $\{t_0,t_1,\ldots\},$ model $\{A_j(x)\}_{j=1}^n, B, R,Q$, $\wt{\mu}^0$, $\Sigma_0$
		\For{$i=0,1,\ldots$}
			\State Solve \eqref{eq: sdre1} and obtain $\Pi(x(t_i);\wt{\mu}^i)$
			\State Set $K(x(t_i);\wt{\mu}^i)$ using \eqref{eq: sdre1}
			\State Set $u(t;\wt{\mu}^i):=-K(x(t_i);\wt{\mu}^i)x(t) \quad \mbox{for }t \in [t_i, t_{i+1}]$ 
			\State Observe the trajectories $x^{i+1}(u(t_i;\wt \mu^i),\mu^*)$
                \State Compute $X$ and $Y$ as in \eqref{eq: XY_full}
                \If{$\|\wt{\mu}^{i}-\wt{\mu}^{i-1}\|_\infty<tol_\mu$}
			      \State Compute $\wt{\mu}^{i+1}$ using BLR
                \Else
                    \State $\wt{\mu}^{i+1}=\wt{\mu}^{i}$
                \EndIf               	
		\EndFor
	\end{algorithmic}
\end{algorithm}

\section{Numerical experiments}\label{sec:test}

We present two numerical test cases with nonlinear PDEs to validate our proposed approach. 
The first test is a nonlinear diffusion-reaction equation, known as the Allen-Cahn equation already presented in Section \ref{sec:motivational}. The second test studies an advection problem with a nonlinear source.

The goal of all our tests is the stabilization of the (unknown) dynamics to the origin by means of the minimization of the given cost functional \aletext{\eqref{disc_costmu}}. Our examples consider problems that do not reach the desired state without the control input.


In both test cases the nonlinear PDE spatial discretization can be written in the form \eqref{eq:dynmu}, and therefore Alg. \ref{alg:RL} and Alg. \ref{alg:RL POD red} can be applied. We will present a comparison of the results and of the computational costs, showing how the POD-DEIM approach can reach results as accurate as the ones provided by Alg. \ref{alg:RL}, but in much less time. The speedup factors might be even more important if we had larger $d$. The number of POD basis and DEIM points will be set as the rank of the provided matrices. In this case, we observed better results towards our aim. Also the classical POD energy in this context is not meaningful since we do not compute optimal snapshots.

We start with an initial parameter configuration $\wt\mu^0=[1,1,1]$ and with a covariance matrix $\Sigma_0 = c\textbf{1}_3$ with $c=10^6$. Since we do not have a priori knowledge on the quality of the initial estimate $\wt\mu^0$, we consider $c$ large (see e.g. \cite{Rasmussen2006}). In Algorithm \ref{alg:RL POD red}, we compute the snapshots set $\mathcal{S}_{lin}$ and, since the parameters are unknown, we apply the strategy described in Section \ref{sec:motivational} for the case of unknown models, and take snapshots from different parameter combinations. We then reduce the model using the POD-DEIM approach and start with the online stage. We recall that we assume trajectories to be observable. In the numerical tests they are integrated in time using an implicit Euler method, accelerated through a Jacobian--Free Newton Krylov method. This must be considered as a ``black box" that, given an initial state, a control, and a time interval $[t_i, t_{i+1}]$, returns the state at time $t_{i+1}$.

Furthermore, in the parameter approximation step we can add noise to the matrix $X$ (or $X_r$) in order to simulate an error on the observation and/or on the model approximation. 
The noise, when added, will be done as follows: for each column $j$ of the matrix $X$ we compute the mean $m_j$ of the absolute values of its elements and choose a positive value $\hat\sigma\in\R$. Then, we add to each component of the $j$-th column of $X$ a gaussian noise with mean $0$ and standard deviation equal to $\hat \sigma m_j$. If $\hat\sigma=0.03$ we will say that noise is $3\%$.


\subsection{Test 1: Allen-Cahn 2D}
In our first test, we consider a PDE in the form \eqref{ex1} where the parameters $\mu_1$, $\mu_2$ and $\mu_3$ are unknown, with the associated cost functional \eqref{costex1}. The true values of the parameters $\mu^*$ are as in Section \ref{sec:motivational}. We apply Alg. \ref{alg:RL} and Alg. \ref{alg:RL POD red} to the problem. As for Section \ref{sec:motivational}, we use $\Delta t =0.025$ and a space grid of $n_{\xi_1} n_{\xi_2}$ nodes with $n_{\xi_1} = n_{\xi_2} = 101$ ($d=10201$). For the reduced model, we computed the snapshots using $\mathcal{S}_{lin}$ described in Section \ref{sec:motivational} with dimensions $r=19$ and $\ell=20$. Our assumption on the parameter space is $\mathcal{D}=[0.1,1]\times[0,15]\times[0,15].$

\begin{figure}[htbp]
    \centering \includegraphics[scale=0.21]{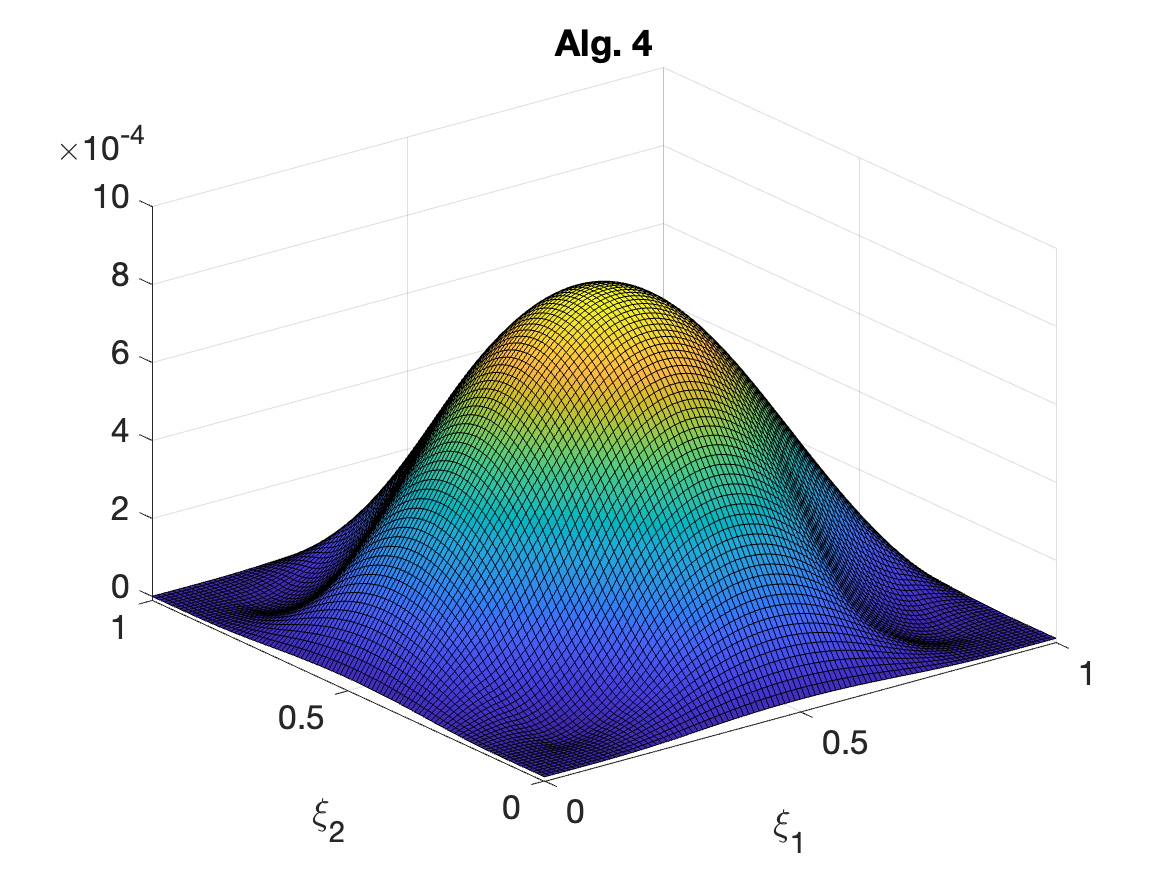}
    \includegraphics[scale=0.21]{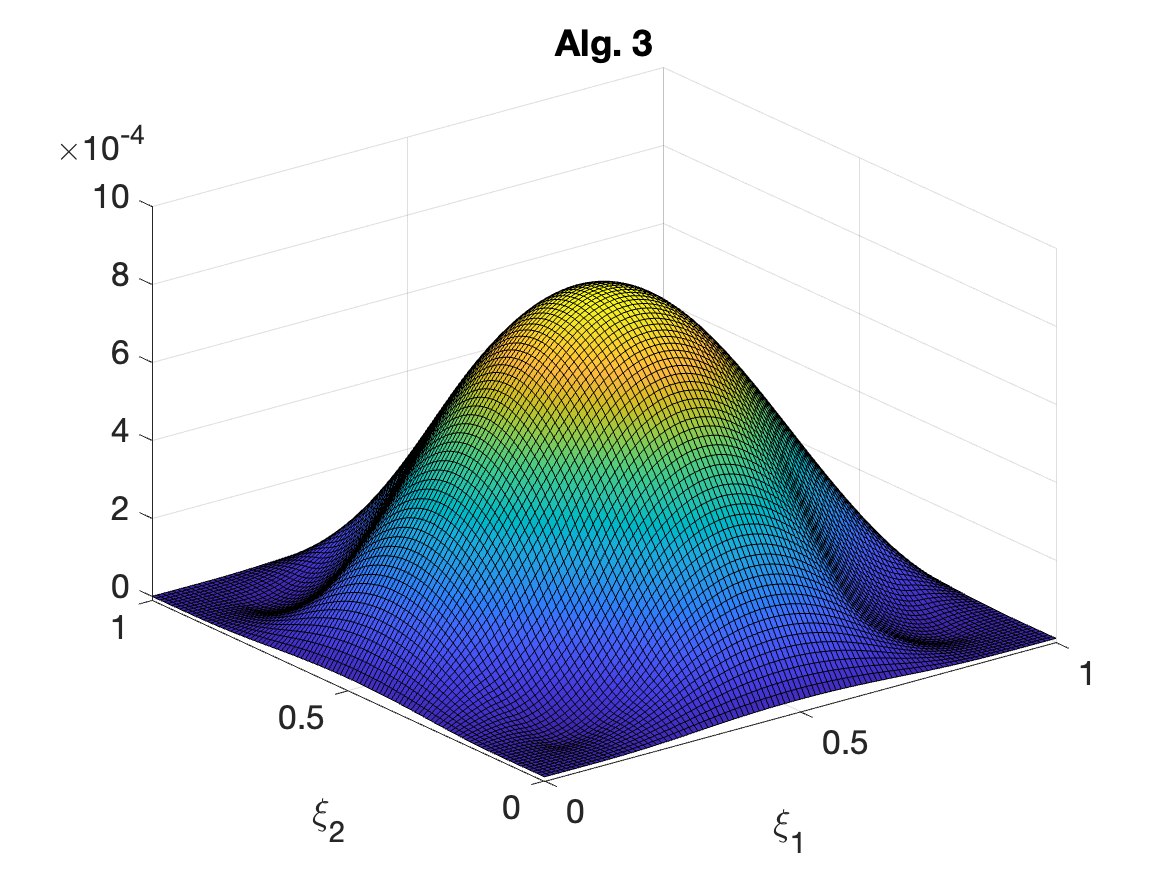}
\includegraphics[scale=0.21]{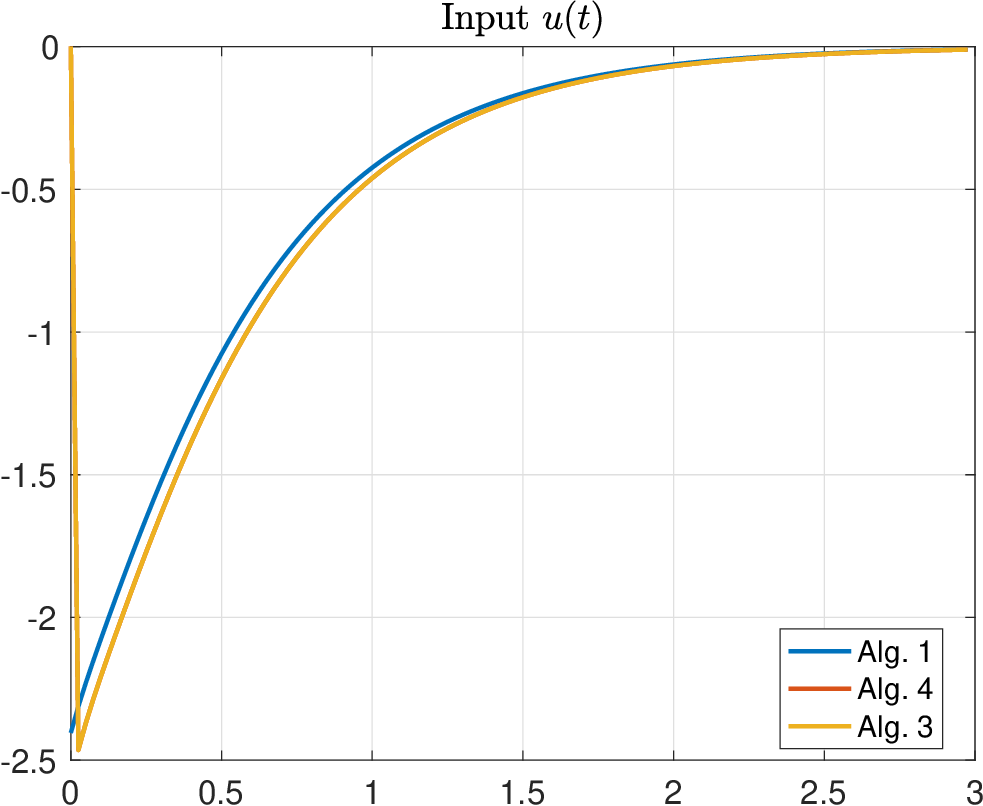}
    \caption{Test 1. Solution at time $t=3$ from Alg. \ref{alg:RL} (left) and Alg. \ref{alg:RL POD red} (middle). See the right panel of Figure \ref{fig:ex1} for a comparison with algorithm \ref{alg: sdre}. The right panel of this picture shows the controls found by the three mentioned algorithms. No noise was added to the data in this test.}
    \label{fig:test1_no_noise}
\end{figure}

Figure \ref{fig:test1_no_noise} shows the results without noise for both Alg. \ref{alg:RL} and Alg. \ref{alg:RL POD red}. In both cases the algorithms stabilize the system and, as shown in the right panel of Figure \ref{fig:test1_no_noise} the controls applied are very similar.


Table \ref{tab:test1_no_noise}  presents the parameter approximation at time $t=3$, without noise. We can observe that the algorithm  using POD-DEIM provides an accurate parameter approximation, even though it starts very far from the desired configuration $\mu^*.$ 


\begin{table}[htbp]
\begin{tabular}{cccccc}
& True $\mu^*$ & $\wt\mu$ & $\wt\mu_r$ & Initial guess $\wt\mu^0$ \\ 
\hline 
$\Delta y$ & \phantom{xx}0.5 & \phantom{xx}0.5 & \phantom{xx}0.5 & 1 \\ 
$y$ & \phantom{x}11 & \phantom{x}11 & \phantom{x}11 & 1 \\ 
$y^3$ & -11 & -10.9998 & -10.9997 & 1 \\ 
\hline 
\end{tabular}
\caption{Test 1: from the left, values of true parameter $\mu^*$, parameter approximation $\wt \mu$ with Alg. \ref{alg:RL} at time $t=3$, parameter approximation $\wt \mu_r$ with Alg. \ref{alg:RL POD red} at time $t=3$ and initial guess corresponding to each term of \eqref{ex1}. No noise was added to the data.}
\label{tab:test1_no_noise}
\end{table}

The benefits of Alg. \ref{alg:RL POD red} are shown in Table \ref{tab:test1_CPUtime} where we compare the CPU times.  Our reference is the CPU needed to execute Alg. \ref{alg: sdre}, that needs a full knowledge of the $\mu^*$ coefficients. It is evident how the POD-DEIM technique used in Alg. \ref{alg:RL POD red} drastically reduced the computational cost of the method with a speed up factor of $60\times$. The time taken to execute Alg. \ref{alg: sdre} is slightly higher than the one of Alg. \ref{alg:RL}. This is because the second one, even if it also performs the parameter approximation, solves one ARE less than \ref{alg:RL}. The time reported in Table \ref{tab:test1_CPUtime} for Alg. \ref{alg:RL POD red} considers only the online stage.


\begin{table}
\begin{tabular}{cccc}
& Algorithm \ref{alg: sdre} & Algorithm \ref{alg:RL} & Algorithm \ref{alg:RL POD red}\\ 
\hline 
CPU time & $130.22$s & $129.59$s & $2.26$s\\ 
\hline 
\end{tabular}
\caption{Test 1. CPU time in seconds of three different algorithms.}
\label{tab:test1_CPUtime}
\end{table}

Finally, in Table \ref{tab:test1_cost}, we provide the evaluation of the cost functionals obtained using different methods. The evaluation of the costs obtained using Alg. \ref{alg:RL POD red} and Alg. \ref{alg:RL} are identic and slightly bigger  than the cost obtained with Algorithm \ref{alg:RL}. This result was expected since SDRE assumes the knowledge of the model and computes directly an input related to the problem.

\begin{table}[htbp]
\begin{tabular}{ccccc}
& Uncontrolled & Algorithm \ref{alg: sdre} & Algorithm \ref{alg:RL} & Algorithm \ref{alg:RL POD red}\\ 
\hline 
Cost Functional & 0.82285 & 0.10924 & 0.11724 & 0.11724\\ 
\hline 
\end{tabular}
\caption{Test 1: evaluation of the cost functional of four different approaches, until time $t=3$. From the left, the first one is the cost of the uncontrolled problem, i.e. the cost of $u(t)\equiv 0$. The second one is the cost obtained using the SDRE method, that uses knowledge of $\mu^*$. This algorithm can return a suboptimal control, but we still use its result as a reference for the optimal cost. The last columns contain the cost of the controls $u$ and $u_r$ found by Alg. \ref{alg:RL} and Alg. \ref{alg:RL POD red}. No noise has been added to this simulation.}
\label{tab:test1_cost}
\end{table}

\begin{figure}[htbp]
    \centering 
    \includegraphics[scale=0.21]{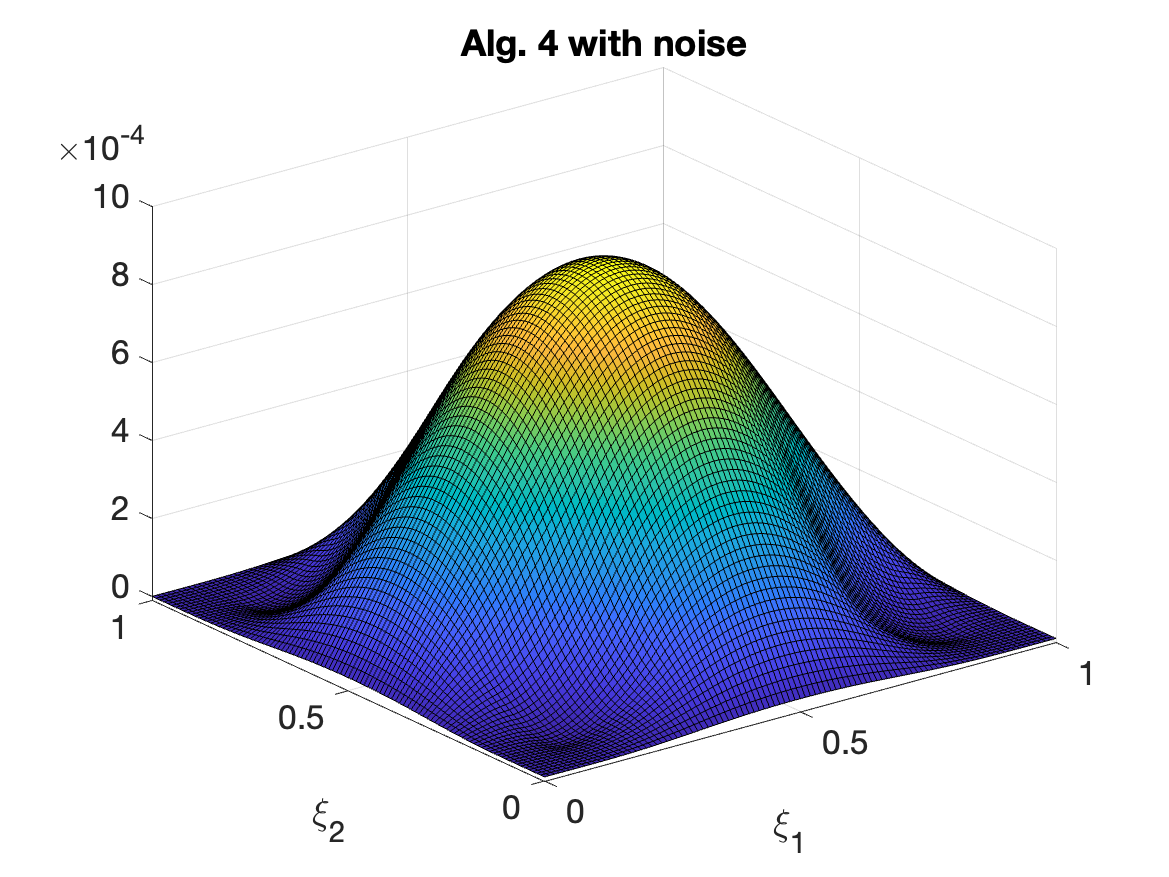}   
\includegraphics[scale=0.21]{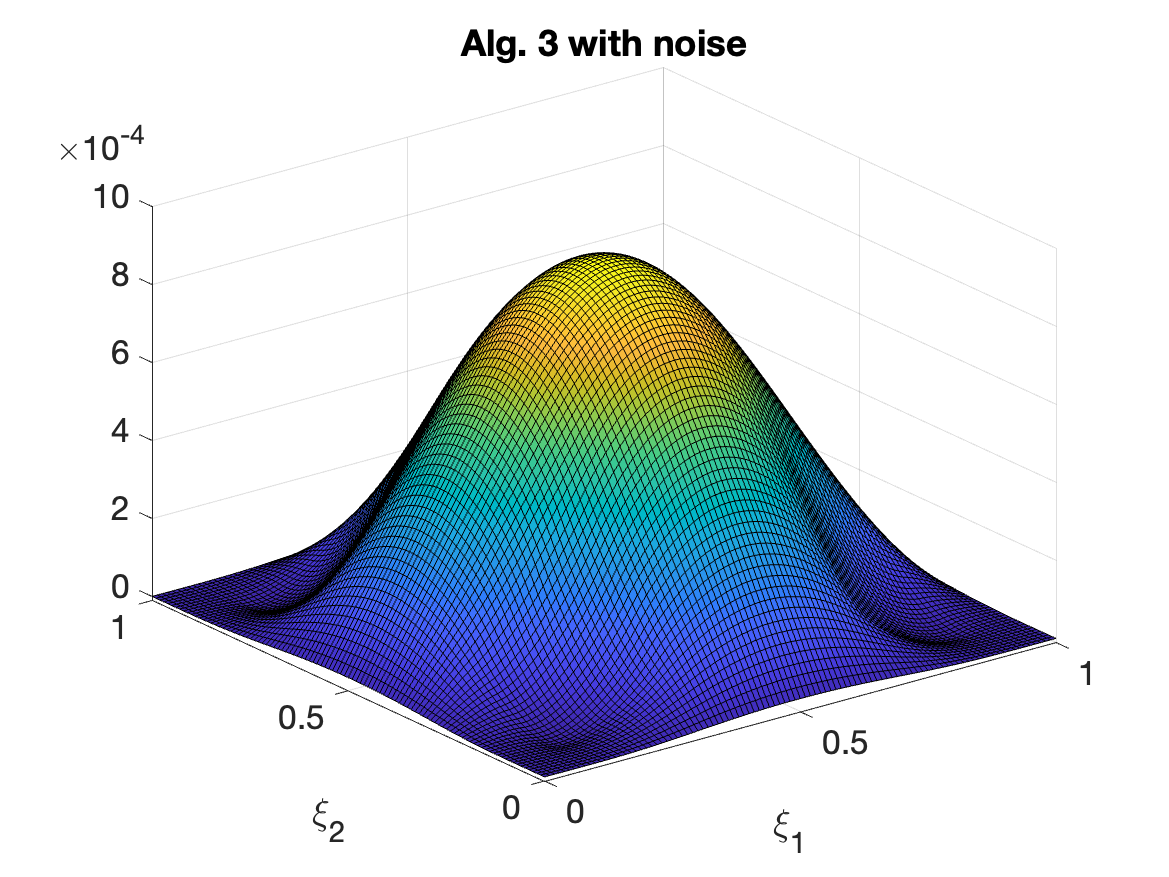}   
\includegraphics[scale=0.21]{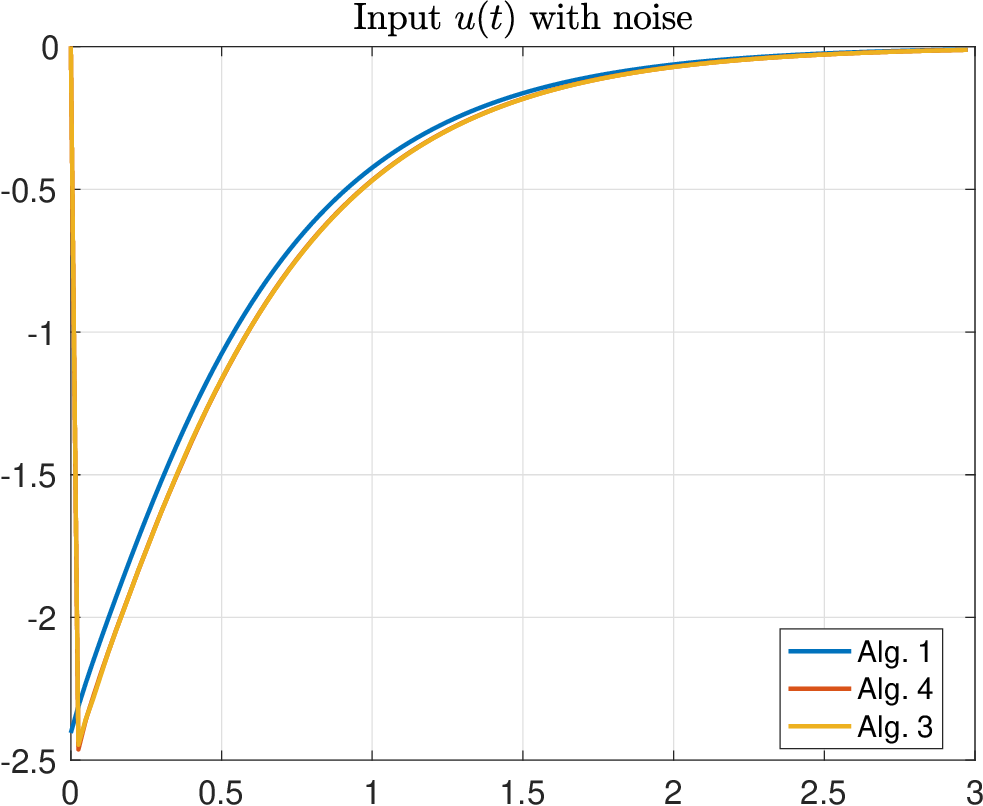} 
    \caption{Test 1. Solution at time $t=3$ from Alg. \ref{alg:RL} (left) and Alg. \ref{alg:RL POD red} (middle) and comparison between the controls found by Alg. \ref{alg: sdre}, Alg. \ref{alg:RL} and  Alg. \ref{alg:RL POD red} (right). Noise of $3\%$ was added to the data.}
    \label{fig:test1_3_noise}
\end{figure}

{\bf Simulations with noise.} To conclude this test we provide the results of our methods when a $3\%$ is added at each iteration. In Figure \ref{fig:test1_3_noise}, we can observe that both Alg \ref{alg:RL POD red} and Alg. \ref{alg:RL} stabilize the problem using similar control inputs.

The parameter estimates are presented in Table \ref{tab:test1_3_noise}. We can observe that the influence of the noise does not allow to find a parameter estimate as good as the one in the results shown in Table \ref{tab:test1_no_noise}. However, it is clear that our proposed method is not far from the desired configuration.

\begin{table}[htbp]
\begin{tabular}{cccccc}
& True $\mu^*$ & $\wt\mu$ & $\wt\mu_r$ & Initial guess $\wt\mu^0$ \\ 
\hline 
$\Delta y$ & 0.5 & \phantom{xx}0.4999 & \phantom{xx}0.4999 & 1 \\ 
$y$ & 11 & \phantom{x}10.9959 & \phantom{x}10.9962 & 1 \\ 
$y^3$ & -11 & -11.0647 & -10.9769 & 1 \\ 
\hline 
\end{tabular}
\caption{Test 1: from the left, values of true parameter $\mu^*$, parameter approximation $\wt \mu$ with Alg. \ref{alg:RL} at time $t=3$, parameter approximation $\wt \mu_r$ with Alg. \ref{alg:RL POD red} at time $t=3$ and initial guess corresponding to each term of \eqref{ex1}. $3\%$ of noise was added to the data.}
\label{tab:test1_3_noise}
\end{table}

\subsection{Test 2: 2D Nonlinear advection equation}

We consider an advection problem with a nonlinear source and zero Dirichlet boundary conditions. We study the following PDE for $ t \in (0, 3)$  and $\xi \in \Omega =(0,1)^2\subset\R^2$ 
\begin{equation}\label{ex2}
    \left\{ \begin{aligned}
    y_t(t,\xi) &= \mu_1 \Delta y(t,\xi) + \mu_2y(t,\xi)\bigl(y_{\xi_1}\!(t,\xi)\! +\! y_{\xi_2}\!(t,\xi)\bigr) + \mu_3y(t,\xi)e^{-0.1y(t,\xi)} + \chi_{\Omega_c}(\xi)u(t)\\
    y(t,\xi) &= 0, \quad \xi\in\partial\Omega\\
    y(0,\xi) &=y_0(\xi), \quad \xi\in\Omega
    \end{aligned} \right.
\end{equation}
where $y_t$, $y_{\xi_1}$ and $y_{\xi_2}$ denote partial derivatives, $y_0(\xi)=0.2sin(\pi\xi_1)sin(\pi\xi_2)$ and $\Omega_c$ is defined as in Section \ref{sec:motivational}. The cost functional is given by \eqref{costex1} and the parameters $(\mu_1, \mu_2, \mu_3)\in\mathcal{D}=[0.1,1]\times[0,1]\times[0,6]$ are unknown. The semi-discretization of this problem can be written in the form \eqref{eq:dynmu} as follows:

\begin{equation}\label{eq:disc2}
\dot x(t) = \Bigl(\mu_1\Delta_d+\mu_2T(x(t))+\mu_3 \text{diag}(e^{-0.1x(t)})\Bigr)x(t) +B u(t)
\end{equation}
where $\Delta_d$, $\text{diag}(\cdot)$ and $B$ are defined as in Section \ref{sec:motivational} and $T(x)$ is a matrix function that approximates the sum of partial derivatives using an upwind scheme.
We note that \eqref{eq:disc2} is in the desired form and so we can study Alg. \ref{alg:RL} and Alg. \ref{alg:RL POD red}.
We use $\Delta t=0.025$ and $n_{\xi_1}=n_{\xi_2}=101$ $(d=10201)$.

\begin{figure}[htbp]
\centering
\includegraphics[scale=0.3]{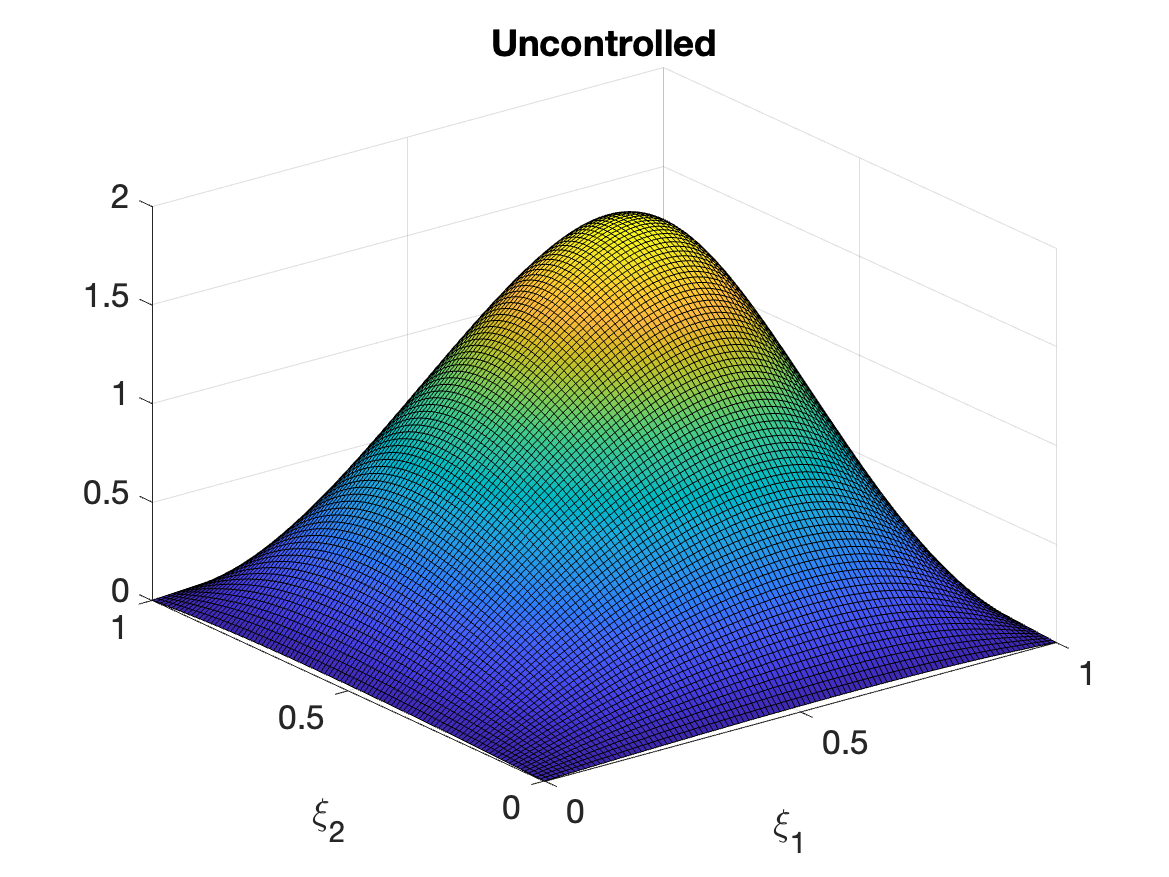}
\includegraphics[scale=0.3]{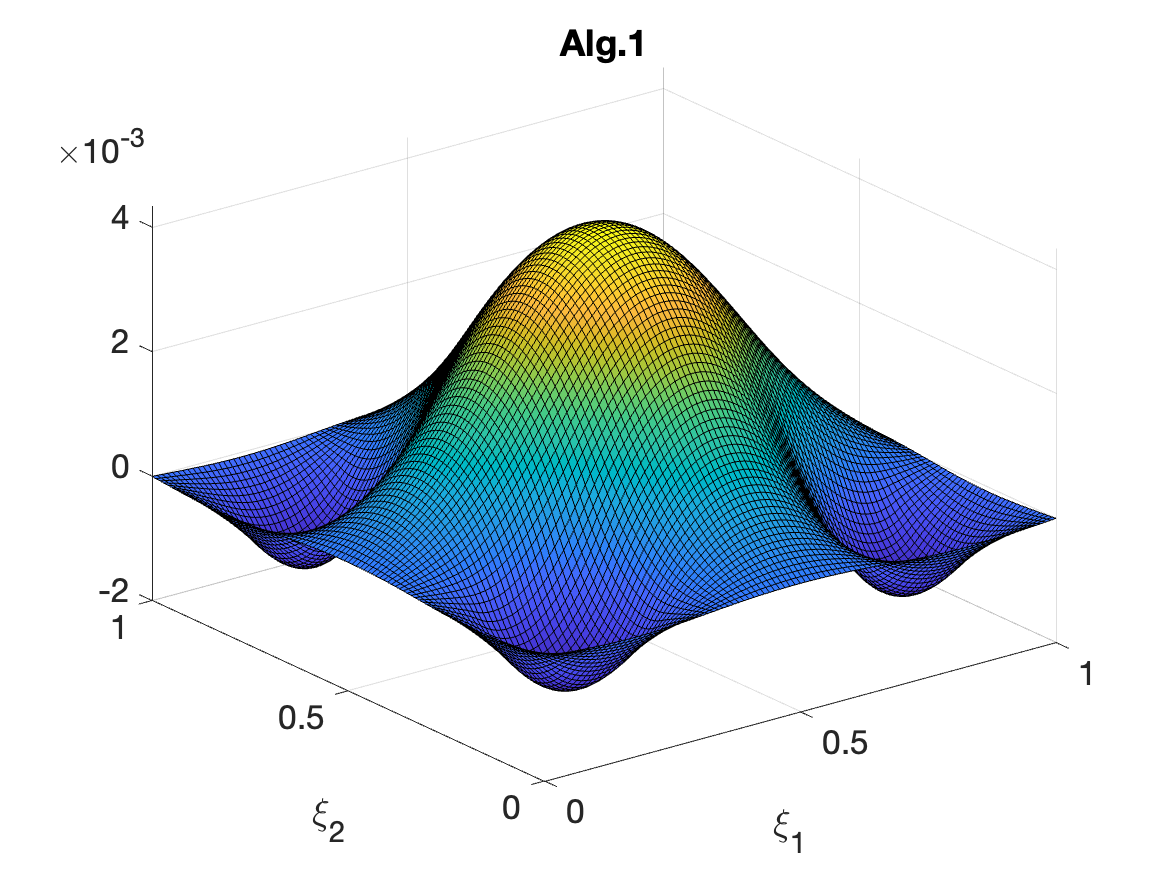}
\caption{Solutions at time $t=3$ of \eqref{ex2} with $\mu_1=0.2,\mu_2=1,\mu_3=5.5$ for the uncontrolled problem $u(t)=0$ (left), and the stabilized solution through Algorithm \ref{alg: sdre} (right). }
\label{fig:ex2}
\end{figure}

We first provide the simulations for the uncontrolled problem in the left panel of Figure \ref{fig:ex2} and the controlled solution using Alg. \ref{alg: sdre} in the right panel. It is clear that Alg. \ref{alg: sdre} stabilizes the solution.

Next, we want to identify and control problem \eqref{ex2} using the following (unknown) parameter configuration $\mu^*_1=0.2,\mu^*_2=1,\mu^*_3=5.5$.
The snapshots are computed using the third approach in Section \ref{sec: sdre}. The linearized equation of \eqref{ex2} around $y(t) \equiv 0$ is:

\begin{equation}\label{ex2_lin}
    \left\{ \begin{aligned}
    y_t(t,\xi) &= \mu_1 \Delta y(t,\xi) + \mu_3y(t,\xi) + \chi_{\Omega_c}(\xi)u(t)\\
    y(t,\xi) &= 0, \quad \xi\in\partial\Omega\\
    y(0,\xi) &=y_0(\xi), \quad \xi\in\Omega
    \end{aligned} \right.
\end{equation}
The finite set of parameters for $\mu_1$ is $I_1=\{0.1,0.2154,0.4642 \}$ whereas for $\mu_3$ we consider $I_3=\{0,1,2,3,4,5,6\}$.

\begin{figure}[htbp]
    \centering \includegraphics[scale=0.21]{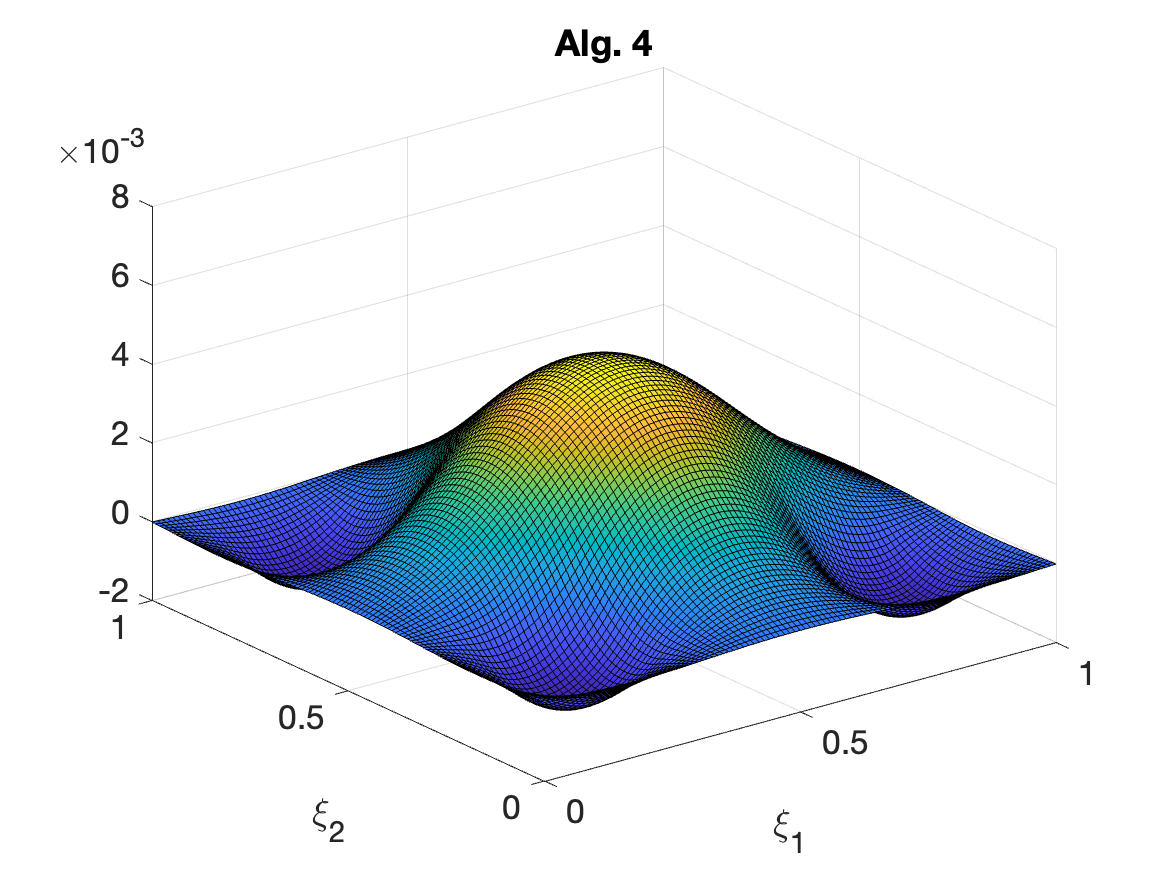}
    \includegraphics[scale=0.21]{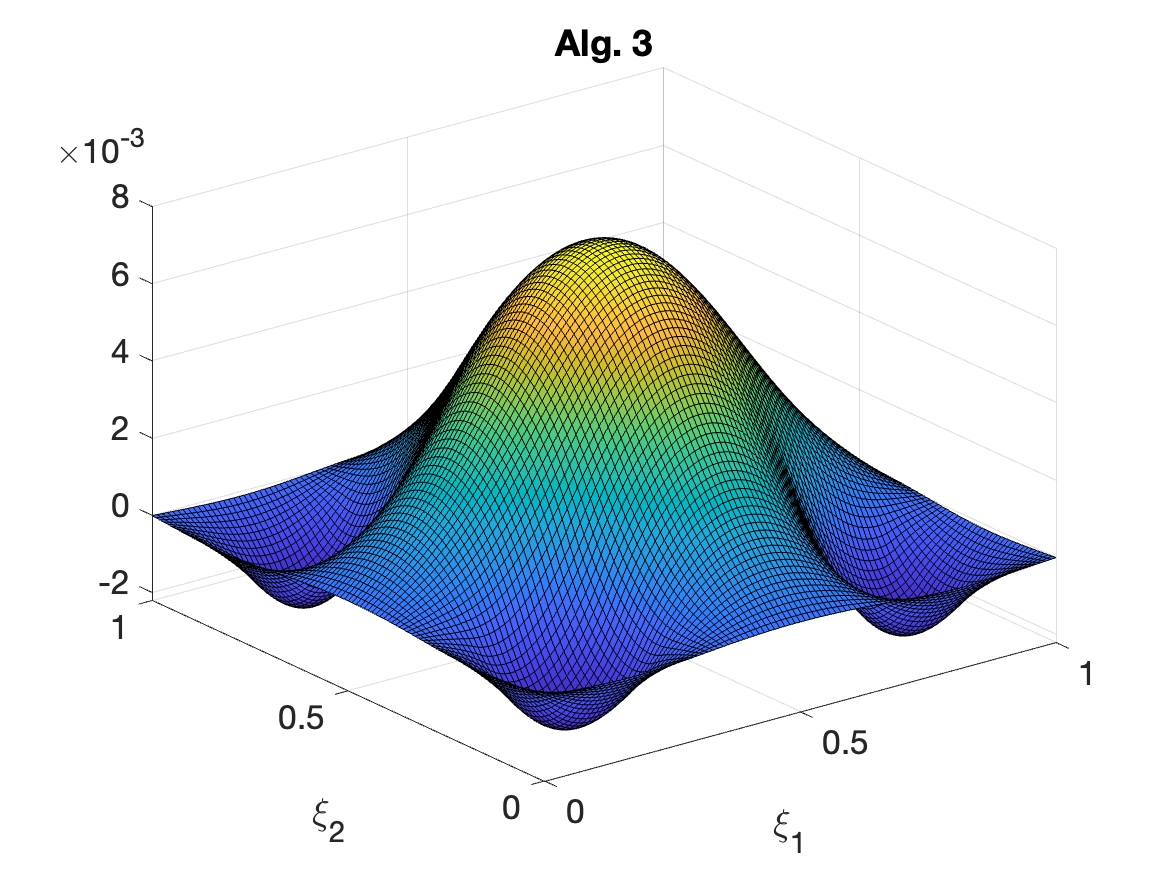}
    \includegraphics[scale=0.21]{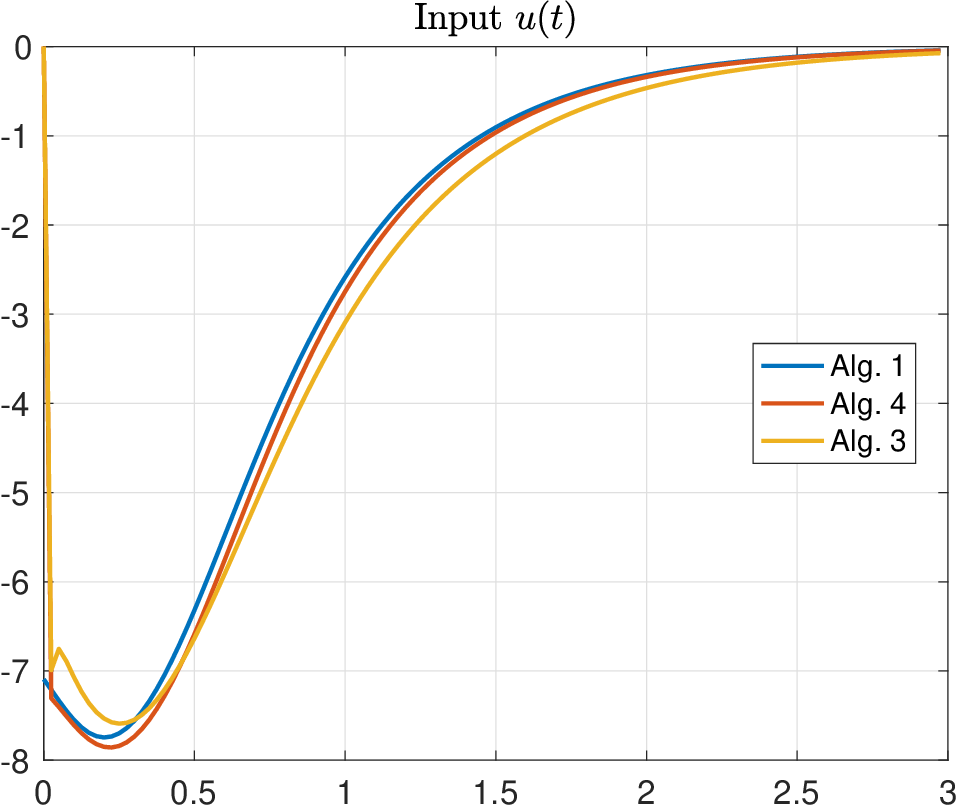}
     \caption{Test 2. Solution at time $t=3$ from Algorithm \ref{alg:RL} (left) and Algorithm \ref{alg:RL POD red} (middle). The right panel shows the controls found by the two mentioned algorithms and Algorithm \ref{alg: sdre}. No noise was added to the data in this test.}
    \label{fig:test2_no_noise}
\end{figure}

We reduce the model with $r=84$ and $\ell=63$. We note that the values of $r$ and $\ell$ are larger in this example with respect to the previous test. This is expected due to the presence of the advection term.
Figure \ref{fig:test2_no_noise}  shows the results of the proposed method, comparing the results of Alg. \ref{alg:RL POD red} and Alg. \ref{alg:RL}. We can observe that we are able to stabilize the problem and that the control inputs are in a very good agreement.

In Table \ref{tab:test2_no_noise}, we provide the parameter configuration estimated from Alg. \ref{alg:RL} and Alg. \ref{alg:RL POD red}. We can observe that Algorithm \ref{alg:RL} (third column) provides the exact desired configuration whereas Alg. \ref{alg:RL POD red} is very close. This is expected since we are adding an error using the surrogate model.


\begin{table}[htbp]
\begin{tabular}{ccccc}
& True $\mu^*$ & $\wt\mu$ & $\wt\mu_r$ & Initial guess $\wt\mu^0$ \\ 
\hline 
$\Delta y$ & 0.2 & 0.2 & 0.2012 & 1 \\ 
$y(y_{\xi_1}+y_{\xi_2})$ & 1 & 1 & 0.9972 & 1 \\ 
$yexp(-0.1y)$ & 5.5 & 5.5 & 5.433 & 1 \\ 
\hline 
\end{tabular}
\caption{Test 2: from the left, values of true parameter $\mu^*$, parameter approximation $\wt \mu$ with Alg. \ref{alg:RL} at time $t=3$, parameter approximation $\wt \mu_r$ with Alg. \ref{alg:RL POD red} at time $t=3$ and initial guess corresponding to each term of \eqref{ex1}. No noise was added to the data.}
\label{tab:test2_no_noise}
\end{table}

In Table \ref{tab:test2_CPUtime}, we provide the computational costs to obtain the results. We can observe that Alg. \ref{alg: sdre} and Alg. \ref{alg:RL}  have similar performances whereas the POD-DEIM approach provides a speed factor of $4\times$. In this example the speedup is smaller compared to the previous test. Indeed, the number  $r$ is larger here, and it is well-known that advection terms are harder to approximate with a POD method.
\begin{table}[htbp]
\centering
\begin{tabular}{ccccc}
& Algorithm \ref{alg: sdre} & Algorithm \ref{alg:RL} & Algorithm \ref{alg:RL POD red}\\  
\hline 
CPU time & $207.3$s & $205.7$s & $50.8$s\\ 
\hline 
\end{tabular}
\caption{Test 2. CPU time in seconds of three different algorithms.}
\label{tab:test2_CPUtime}
\end{table}

Finally, in Table \ref{tab:test2_cost} we provide the evaluation of the cost functional for the three compared algorithms. As expected, we can observe that the value of the cost using Alg. \ref{alg:RL POD red} is slightly larger than the other methods but very comparable.

\begin{table}[htbp]
\centering
\begin{tabular}{ccccc}
& $u=0$ & Algorithm \ref{alg: sdre} & Algorithm \ref{alg:RL} & Algorithm \ref{alg:RL POD red}\\ 
\hline 
Cost Functional & 10.3886 & 2.6136 & 2.7384 & 2.792\\ 
\hline 
\end{tabular}
\caption{Test 2: evaluation of the  cost funcitonal for four different controls, until time $t=3$. No noise was added.}
\label{tab:test2_cost}
\end{table}

Overall, we can conclude that POD-DEIM method obtains results extremely close to the ones obtained by Alg. \ref{alg:RL}.However, the computational cost is significantly reduced.\\

\begin{figure}[htbp]
    \centering 
    \includegraphics[scale=0.21]{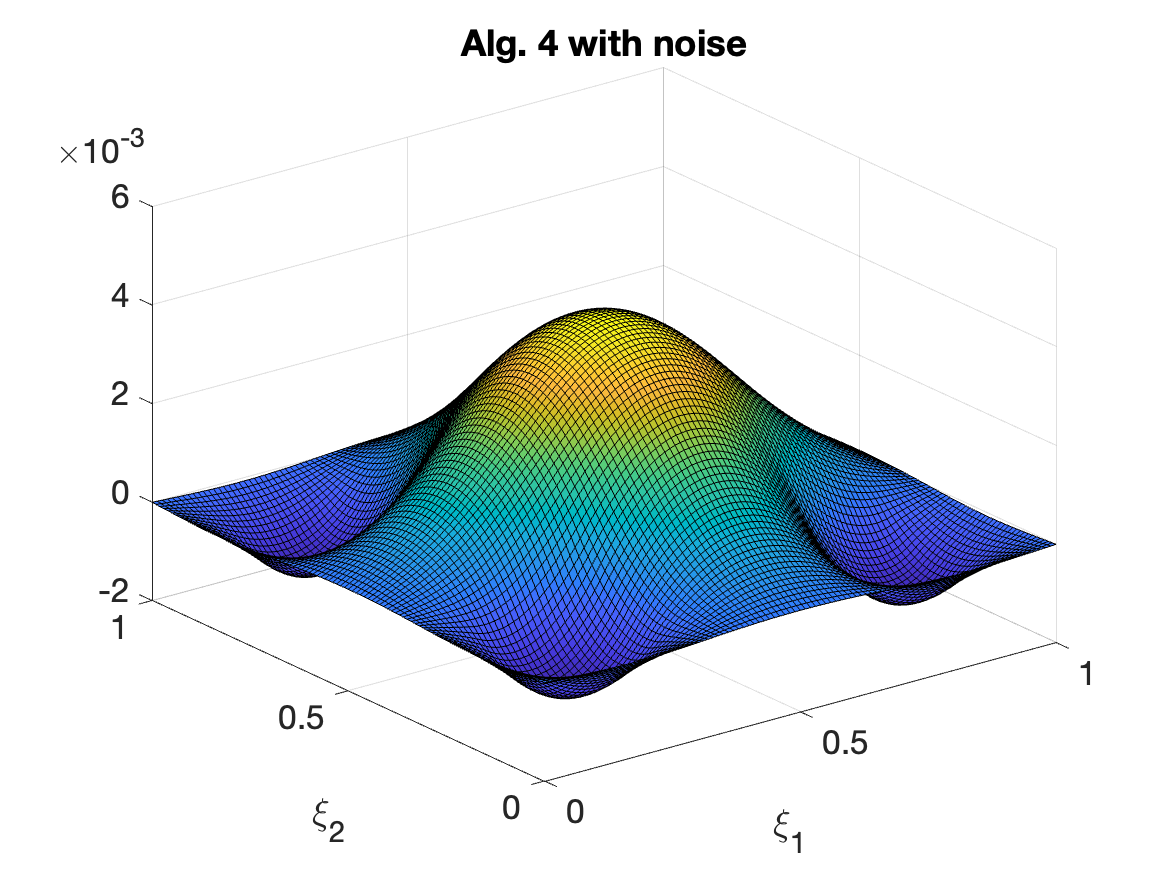}    
\includegraphics[scale=0.21]{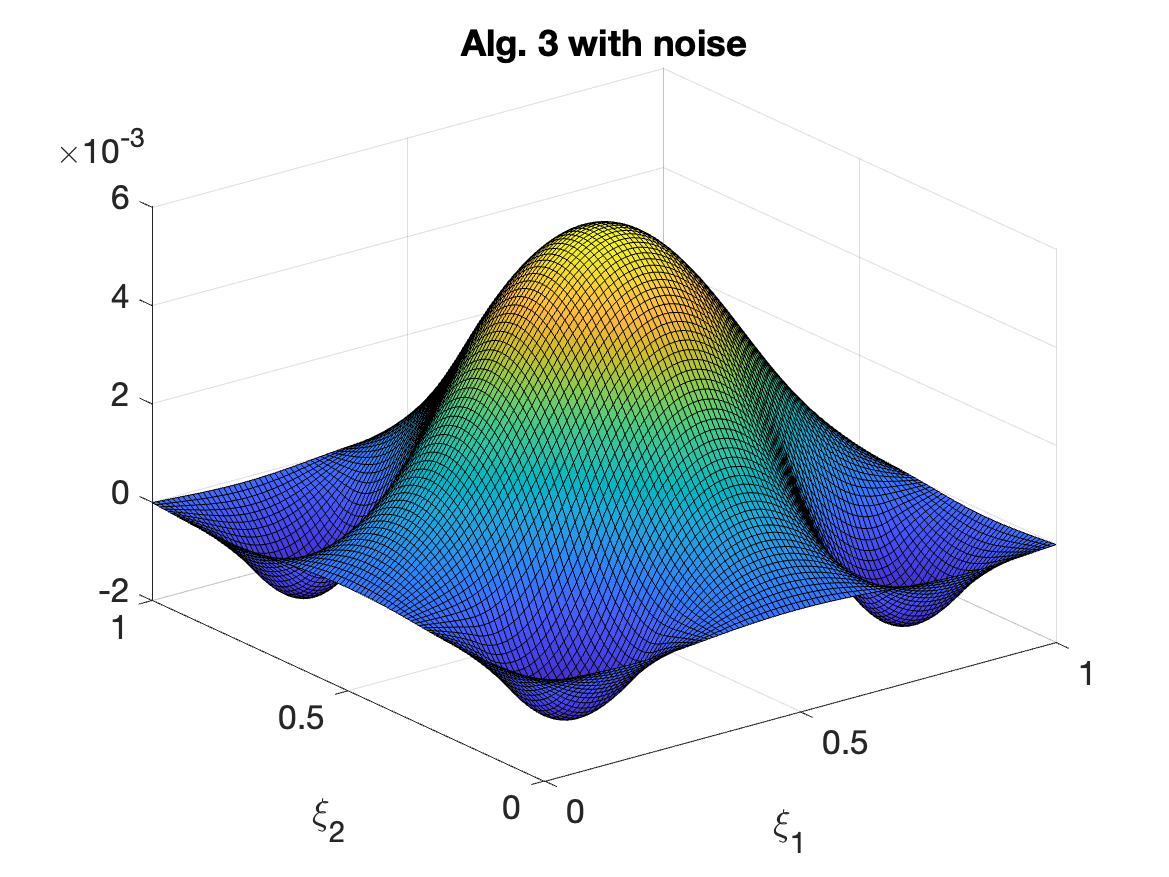}    
 \includegraphics[scale=0.21]{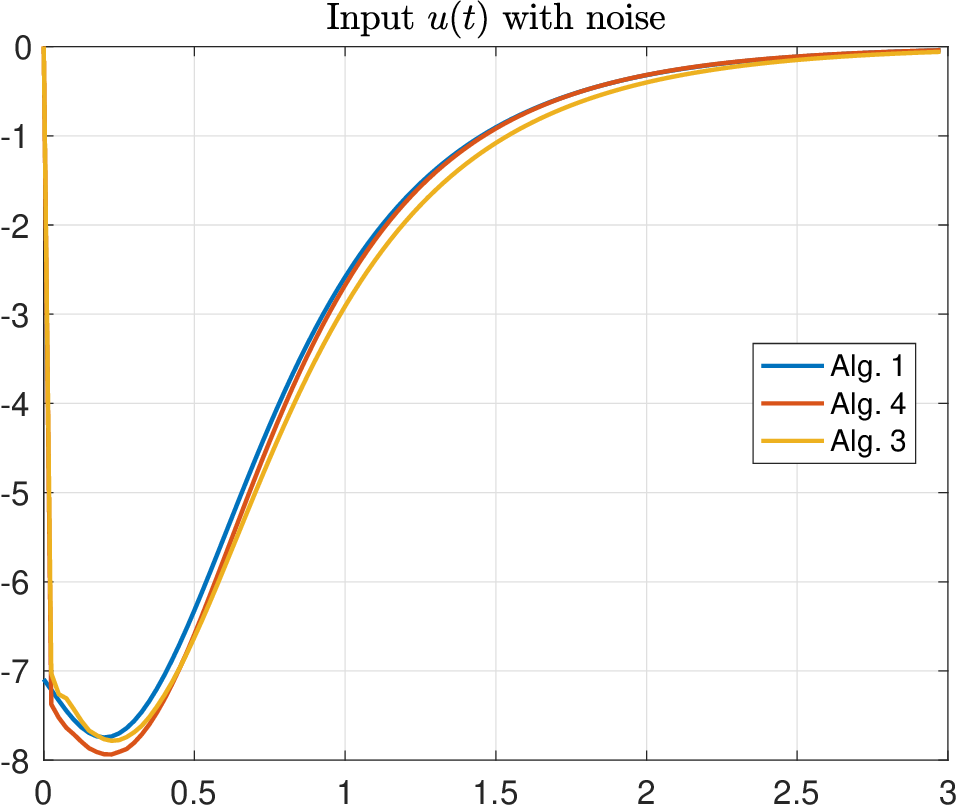}
   \caption{Test 2. Solution at time $t=3$ from Alg. \ref{alg:RL} (left) and Alg. \ref{alg:RL POD red} (middle) and comparison between the controls found by Alg. \ref{alg: sdre}, Alg. \ref{alg:RL} and Alg. \ref{alg:RL POD red} (right). Noise of $3\%$ was added to the data. }
    \label{fig:test2_3_noise}
\end{figure}

{\bf Simulations with noise.} To conclude this test, we provide the results of our methods when a $3\%$ of noise is added each iteration, i.e. with $\hat\sigma=0.03$. In Figure \ref{fig:test2_3_noise}, we can observe that both Alg \ref{alg:RL POD red} and Alg. \ref{alg:RL} stabilize the problem using similar control inputs.

The estimated parameter under the presence of noise is presented in Table \ref{tab:test2_3_noise}. As in the previous test, we can observe that the influence of the noise does not allow to estimate the parameter as accurately as in Table \ref{tab:test2_no_noise}. Nevertheless, the obtained result is not far from the true configuration.


\begin{table}
\begin{tabular}{ccccc}
& True $\mu^*$ & $\wt\mu$ & $\wt\mu_r$ & Initial guess $\wt\mu^0$ \\ 
\hline 
$y_{xx}$ & 0.2 & 0.1996 & 0.2006 & 1 \\ 
$yy_x$ & 1 & 1.0013 & 1.0024 & 1 \\ 
$yexp(-.1y)$ & 5.5 & 5.5149 & 5.4664 & 1 \\ 
\hline 
\end{tabular}
\caption{Test 2: from the left, values of true parameter $\mu^*$, parameter approximation $\wt \mu$ with Alg. \ref{alg:RL} at time $t=3$, parameter approximation $\wt \mu_r$ with Alg. \ref{alg:RL POD red} at time $t=3$ and initial guess corresponding to each term of \eqref{ex1}. $3\%$ of noise was added to the data.}
\label{tab:test2_3_noise}
\end{table}

\section{Conclusions}
\label{sec:conclusions}

We have proposed a new online algorithm that, efficiently, stabilizes and identifies high dimensional problems. The system we want to control is assumed to be unknown but observable through its trajectories. Specifically, the system is considered fully identified based on its parameter configuration.
The algorithm starts with a prior distribution over the unknown parameter and iteratively estimates the parameter configuration, using BLR, thanks to the possibility to observe the dynamics evolution for a given control input. At the same time, we aim at computing a control that minimizes the given cost functional. The control is computed at each time step using the SDRE approach with the current parameter estimate. 

Furthermore, our algorithm reduces the computational cost of the method by using a POD-DEIM approach. The discretization of PDEs in two dimensions leads to high dimensional systems, and the algorithm requires the solution of many ARE. 
 POD-DEIM technique allows to obtain a surrogate model 
 that accelerates the computations. We have also introduced the use of POD-DEIM for SDRE with a keen focus on the choice of the snapshots set which is crucial to obtain suitable basis for the projection.
 Numerical experiments show that the use of POD-DEIM in the identification and control algorithm stabilizes the system and identifies the model, even when we consider noisy observations. Furthermore, we found impressive speedup factors when using model reduction. All the numerical tests exhibit the stabilization of the considered problem and numerical convergence of the parameter estimate towards the true parameters' values. The proof of convergence of the utilized method under suitable assumptions is still an open problem and we plan to study it in the future.

\section*{Declarations}

\bmhead{ Funding} 
A. Alla and A. Pacifico are members of the INdAM-GNCS activity group. The work of A.A. has been carried out within the “Data-driven discovery and control of multi-scale interacting artificial agent systems”, and received funding from the European Union Next-GenerationEU - National Recovery and Resilience Plan (NRRP) – MISSION 4 COMPONENT 2, INVESTIMENT 1.1 Fondo per il Programma Nazionale di Ricerca e Progetti di Rilevante Interesse Nazionale (PRIN) – CUP H53D23008920001. This manuscript reflects only the authors’ views and opinions, neither the European Union nor the European Commission can be  considered responsible for them. The work of A.P. has been supported by MIUR  with  PRIN project 2022 funds (2022238YY5, entitled
"Optimal control problems: analysis, approximation ")


\bmhead{Ethical Approval}
Not applicable.

\bmhead{Authors' contributions}
A.A. designed and directed the project.
A.A and A.P. wrote the paper, developed the numerical tests, created figures and tables presented in the manuscript.

\bmhead{Availability of data and materials}  No Data associated in the manuscript.

\bmhead{ Conflict of interest} No conflict of interest.

\backmatter









\bibliography{bibsdre.bib,Pesare_PhD_Thesis.bib,PODbib.bib}

\end{document}